\documentclass[pdflatex,sn-mathphys-num]{sn-jnl}% Math and Physical Sciences Numbered Reference Style
\usepackage{natbib}%
\usepackage{graphicx}%
\usepackage{multirow}%
\usepackage{amsmath,amssymb,amsfonts}%
\usepackage{amsthm}%
\usepackage{mathrsfs}%
\usepackage[title]{appendix}%
\usepackage{xcolor}%
\usepackage{textcomp}%
\usepackage{manyfoot}%
\usepackage{booktabs}%
\usepackage{algorithm}%
\usepackage{algorithmicx}%
\usepackage{algpseudocode}%
\usepackage{listings}%
\theoremstyle{thmstyleone}%
\newtheorem{theorem}{Theorem}%  meant for continuous numbers
\newtheorem{proposition}[theorem]{Proposition}% 

\theoremstyle{thmstyletwo}%
\newtheorem{remark}{Remark}%
\newtheorem{lemma}{Lemma}%
\newtheorem{notation}{Notation}%

\theoremstyle{thmstylethree}%
\newtheorem{definition}{Definition}%

\usepackage{tikz}

\usepackage{dsfont}

\usepackage{enumitem} % Nécessaire pour personnaliser les listes

\def \d {\mathrm{d}}
\newcommand{\bs}[1]{\boldsymbol{#1}}

\begin{document}

\title[Control of the Fisher-Stefan system]{Control of the Fisher-Stefan system}


\author[1]{\fnm{Idriss } \sur{Boutaayamou}}\email{d.boutaayamou@uiz.ac.ma}

\author[2]{\fnm{Fouad} \sur{Et-tahri}}\email{fouad.et-tahri@edu.uiz.ac.ma}

\author[3,4]{\fnm{Lahcen} \sur{Maniar}}\email{maniar@uca.ma}

\author*[5]{\fnm{Francisco} \sur{Periago}}\email{f.periago@upct.es}

\affil[1]{\orgdiv{Lab-SIV}, \orgname{Polydisciplinary Faculty-Ouarzazate, Ibnou Zohr University}, \orgaddress{
%		\street{Street}, 
\city{Ouarzazate}, \postcode{45000}, \state{B.P. 638}, \country{Morocco}}}

\affil[2]{\orgdiv{Lab-SIV}, \orgname{Faculty of Sciences-Agadir, Ibnou Zohr University}, \orgaddress{
%		\street{Street}, 
\city{Agadir}, 
%\postcode{10587}, 
\state{B.P. 8106}, \country{Morocco}}}

\affil[3]{\orgdiv{University Cadi Ayyad}, \orgname{Faculty of Sciences Semlalia, Laboratory of Mathematics Modelisation and Automatic Systems}, \orgaddress{
		%		\street{Street}, 
		\city{Marrakesh}, \postcode{40000}, \state{B.P. 2390}, \country{Morocco}}}

\affil[4]{\orgdiv{The UM6P Vanguard Center Mohammed VI Polytechnic University}, 
%\orgname{Organization}, 
\orgaddress{\street{Hay Moulay Rachid},
		 \city{Ben Guerir}, 
		 \postcode{43150}, 
%		 \state{State}, 
\country{Morocco}}}

\affil*[5]{\orgdiv{Department of Applied Mathematics and Statistics. Technical University of Cartagena}, \orgname{Organization}, 
	\orgaddress{
%		\street{Street}, 
\city{Campus Muralla del Mar}, \postcode{30202}, \state{Cartagena (Murcia)}, \country{Spain}}}


\abstract{This paper addresses the exact controllability of trajectories in the one-dimensional Fisher-Stefan problem---a reaction-diffusion equation that models the spatial propagation of biological, chemical, or physical populations within a free-end domain, governed by Stefan's law. We establish the local exact controllability to the trajectories by reformulating the problem as the local null controllability of a nonlinear system with distributed controls. Our approach leverages the Lyusternik-Graves theorem to achieve local inversion, leading to the desired controllability result. Finally, we illustrate our theoretical findings through several numerical experiments based on the Physics-Informed Neural Networks (PINNs) approach.}




\maketitle

\section{Introduction and main results}	\label{sec1}
In mathematical biology, the one-dimensional prototype model of the \textbf{Fisher-KPP} (Kolmogorov-Petrovsky-Piskunov) equation is formulated as:
\begin{eqnarray}
	\varrho_{t}-\sigma\varrho_{xx} =r\varrho\left(1-\frac{\varrho}{k}\right), \label{FK} 
\end{eqnarray}
where $\sigma>0$ is the diffusivity coefficient, $r$ is the growth rate, and $k>0$ is the carrying capacity density. The model \eqref{FK}
elucidates the space-time dynamics of a positive population density $\varrho(x,t)$, evolving under the effect of a linear diffusion combined with a Logistic source term. See \cite{ canosa1973nonlinear, grindrod1991patterns, kolmogorov1937study, murray2002mathematical} and the references therein. 

In this paper, we are interested in a more general situation where the dynamics governed by \eqref{FK} takes place in a right-moving boundary domain $(0,L(t))$. We will assume that the moving boundary evolves according to the following law:
\begin{eqnarray}
	L^{\prime}(t)=-\mu\varrho_{x}(L(t),t). \label{Stefan}
\end{eqnarray}
Notably, equation \eqref{Stefan} resembles the classical Stefan condition \cite{crank1984free} used to describe the phase transition of materials.
The parameter $\mu$, which corresponds to the inverse Stefan number, serves as a connection between the density gradient at $L(t)$ and the velocity of the boundary. This type of moving-boundary problems have a long history as they have been widely used to simulate a number of physical and industrial processes \cite{brosa2019extended, dalwadi2020mathematical, mitchell2010improving, mitchell2009finite, planella2021extended}.  In recent years, there has been growing interest in using moving-boundary problems to study biological processes, such as cancer cell invasion, cell movement, and wound healing \cite{fadai2020new, gaffney1999modelling, kimpton2013multiple,  shuttleworth2019multiscale, ward1997mathematical}.
\par 
All in all, the following system, known as the \textbf{Fisher-Stefan} model
\begin{equation} \label{s0}
	\left\{
	\begin{aligned}
		&\varrho_{t}-\sigma\varrho_{xx} =r\varrho\left(1-\frac{\varrho}{k}\right) & & \text {in}\; (0,L(t))\times (0,\infty), \\
		&\varrho_{x}(0,t)=0 & & \text {in}\;(0,\infty), \\
		& \varrho(L(t),t)=0  & & \text {in}\;(0,\infty), \\
		&L^{\prime}(t)=-\mu\varrho_{x}(L(t),t) & & \text {in}\;(0,\infty),  \\	
		& L(0)=L_{0},\\
		& \varrho(\cdot,0)=\varrho_{0} & &\text {in}\;(0,L_{0}),
	\end{aligned}
	\right.
\end{equation}
was first proposed by Du and Lin \cite{du2010spreading}, and Du and Guo \cite{du2011spreading, du2012stefan} . 
%	in case uncontrolled $u=0$ and an infinite horizon $T=\infty$, 
In particular, Du and Lin \cite{du2010spreading} showed that \eqref{s0} has solutions whose population density progresses to a traveling-wave solution  with an asymptotic speed that depends on $\mu$. This corresponds to population survival and successful invasion, with
$\varrho(x,t)\rightarrow 1$ when $t\rightarrow\infty$. They also proved that \eqref{s0} admits solutions where the population does not settle, so that $\varrho(x,t)\rightarrow 0$ when $t\rightarrow\infty$.\\
In the context of microbial species, the goal is to apply a control $u(t)$ such that the microbial density $\rho(x,t)$ vanishes across the domain $x\in (0,L(T))$ at a final time horizon $T>0$. Thus, controlling $L(t)$ ensures the complete removal of the microbe. This has important applications in medicine, agriculture, and environmental management. It also falls within the scope of controlling harmful populations or optimizing disinfection. Consequently, we consider \eqref{s0} over a finite time horizon, namely,
\begin{equation} \label{ss1}
	\left\{
	\begin{aligned}
		&\varrho_{t}-\sigma\varrho_{xx} =r\varrho\left(1-\frac{\varrho}{k}\right) & & \text {in}\; (0,L(t))\times (0,T), \\
		&\varrho_{x}(0,t)=u(t) & & \text {in}\;(0,T), \\
		& \varrho(L(t),t)=0  & & \text {in}\;(0,T), \\
		&L^{\prime}(t)=-\mu\varrho_{x}(L(t),t) & & \text {in}\;(0,T),  \\	
		& L(0)=L_{0},\\
		& \varrho(\cdot,0)=\varrho_{0} & &\text {in}\;(0,L_{0}),
	\end{aligned}
	\right.
\end{equation}
where $(\varrho_0, L_0)$ is the initial state of the system and the control $u(t)$, which is applied at $x=0$, acts as a disinfection flux, impacting the microbial density throughout the domain.
In principle, $u(t)$ can be achieved by influencing the boundary flux through $u(t)$ and utilizing diffusion properties to propagate this effect across the entire domain. This type of control might, for instance, represent medical treatment, where antimicrobial substances are administered at a specific point and diffuse through the affected tissue or medium. In agriculture, it could be used to control pathogenic bacteria or fungi in a field or plantation. 
For more information regarding mathematical models of propagation and control strategies for invasive or harmful species \cite{lewis2016mathematics}, and, as far as biological and environmental applications of optimal control are concerned \cite{trelat2005controle}.
\par
%To the best of our knowledge, the controllability of the \textbf{Fisher-Stefan} model \eqref{ss1} has not been explored, 

Focusing on mathematical controllability properties of free-boundary problem, several recent works \cite{barcena2023exact, demarque2018local, fernandez2017local, fernandez2019local, fernandez2016controllability} have investigated the controllability of one-dimensional parabolic free boundary problems. In \cite{barcena2023exact}, a local exact controllability to the trajectories is proved for a system similar to \eqref{ss1}, with Dirichlet boundary control, but without the nonlinear Logistic term. In \cite{demarque2018local, fernandez2017local, fernandez2019local, fernandez2016controllability}, a local null controllability result for the density is established via a distributed control approach.  In \cite{fernandez2017local, fernandez2016controllability}, the proof is based on fixing the free boundary $L\in C^1([0,T])$ and eliminating its dynamics, as well as establishing an observability inequality for the linear heat equation in the non-cylindrical domain $(0,L(t))\times (0,T)$, with an observability constant uniform in $L$. To recover the boundary law, a Schauder fixed-point argument is applied while preserving this law in an appropriate space. In \cite{fernandez2019local}, the same local null controllability result is obtained using a different approach, which is based on a non-trivial transformation to a fixed cylindrical domain, as well as a local inversion argument, while in \cite{barcena2023exact} the technique is based on trivial transformation to a fixed cylindrical domain. 

Nonetheless, these results do not ensure that the density is zero over the entire spatial domain. Therefore, controlling the component $L$ is also relevant for managing the distribution of the density.

In this paper, we address the problem of controlling both the density and the free boundary to trajectories. Our approach to tackling the problem is similar to that in \cite{barcena2023exact}, but in this case we ought to deal with a nonlinear Logistic term. To this end, we rely on Lyusternik-Graves theorem. 

The numerical approximation of the nonlinear system and its linearized counterpart is addressed as well. In this case, we rely on recent deep-learning-based techniques, specifically Physics-Informed Neural Networks (PINNs).

\subsection{Main result and reformulation of the problem}
Throughout the paper, for the sake of simplicity, we will assume 
$$\sigma=r=k=1.$$
The dynamics of the system becomes
\begin{equation} \label{s1}
	\left\{
	\begin{aligned}
		&\varrho_{t}-\varrho_{xx} =\varrho\left(1-\varrho\right) & & \text {in}\; Q_L, \\
		&\varrho_{x}(0,t)=u(t) & & \text {in}\;(0,T), \\
		& \varrho(L(t),t)=0  & & \text {in}\;(0,T), \\
		&L^{\prime}(t)=-\mu\varrho_{x}(L(t),t) & & \text {in}\;(0,T),  \\	
		& L(0)=L_{0},\\
		& \varrho(\cdot,0)=\varrho_{0} & &\text {in}\;(0,L_{0}),
	\end{aligned}
	\right.
\end{equation}
where $Q_L:=(0,L(t))\times (0,T)$. In the following, we consider the energy space:
\begin{eqnarray*}
	H^{2,1}(Q_{L}):=\left\{\varrho\in L^{2}(Q_{L})\;:\; \varrho_{x},\varrho_{xx},\, \varrho_{t}\in L^{2}(Q_{L})\right\}
\end{eqnarray*}
and we use the definition of smooth trajectories of \eqref{s1}:
\begin{definition}
	We call a smooth trajectory of \eqref{s1}, any triplet $(\overline{\varrho},\overline{L},\overline{u})$ belonging to $$\left[H^{2,1}(Q_{\overline{L}})\cap W^{1,\infty}(0,T;H^{1}(0,\overline{L}(t)))\right]\times W^{1,\infty}(0,T)\times H^{1/4}(0,T),$$
	\textcolor{red}{where, for all $t \in [0,T]$, $\overline{L}(t) > L_\star$ for some fixed constant $L_\star > 0$
	} and satisfies:
	\begin{equation*}
		\left\{
		\begin{aligned}
			&\overline{\varrho}_{t}-\overline{\varrho}_{xx} =\overline{\varrho}\left(1-\overline{\varrho}\right) & & \text {in}\; (0,\overline{L}(t))\times (0,T), \\
			&\overline{\varrho}_{x}(0,t)=\overline{u}(t) & & \text {in}\;(0,T), \\
			& \overline{\varrho}(\overline{L}(t),t)=0  & & \text {in}\;(0,T), \\
			&\overline{L}^{\prime}(t)=-\mu\overline{\varrho}_{x}(\overline{L}(t),t) & & \text {in}\;(0,T),  \\	
			& \overline{L}(0)=\overline{L}_{0},\\
			& \overline{\varrho}(\cdot,0)=\overline{\varrho}_{0} & &\text {in}\;(0,\overline{L}_{0}).
		\end{aligned}
		\right.
	\end{equation*}
\end{definition}
Our main finding is as follows
\begin{theorem} \label{main result}
	Let $(\bar{\varrho}, \bar{L}, \bar{u})$ be a smooth trajectory of \eqref{s1}. Then, there
	exists $\varepsilon> 0$ such that for any $L_{0}>\textcolor{red}{ L_{\star}}$ and any $\varrho_{0}\in H^{1}(0,L_{0})$ with $\varrho_{0}(L_{0})=0$ satisfying
	\begin{eqnarray*}
		|L_{0}-\bar{L}(0)|+ \|\varrho_{0}-\bar{\varrho}(\cdot,0)\|_{H^{1}(0,L_{0})}< \varepsilon,
	\end{eqnarray*}
	there exists a control $u\in H^{1/4}(0,T)$ such that the associated state $(\varrho, L)$ belonging to $H^{2,1}(Q_{L})\times H^{1}(0,T)$ with 
	\begin{eqnarray*}
		u(0)=\varrho_{0,x}(0) \quad\mbox{and}\quad L>\textcolor{red}{ L_{\star}},
	\end{eqnarray*}
	and satisfies: 
	\begin{eqnarray*}
		\varrho(\cdot,T)=\bar{\varrho}(\cdot,T) \;\mbox{in}\; (0,\bar{L}(T)) \quad \mbox{and} \quad L(T)=\bar{L}(T). 
	\end{eqnarray*}
	\qed
\end{theorem}
In our study, we reformulate problem \eqref{s1} as a nonlinear parabolic equation on a fixed cylindrical domain by the change of state:
\begin{eqnarray}
	\psi(x,t):=\varrho(x L(t),t)\;\mbox{and}\; h(t):=L^{2}(t)\; \mbox{for}\;(x,t)\in Q:=(0,1)\times (0,T). \label{Chang}
\end{eqnarray}
By applying the chain rule, we obtain the following system of equations for $(\psi,h)$:
\begin{equation} \label{s2}
	\left\{
	\begin{aligned}
		&h\psi_{t}-\psi_{xx} - \frac{x}{2}h^{\prime}\psi_{x}=h\psi(1-\psi) & & \text {in}\; Q, \\
		& \psi_{x}(0,\cdot)=\textcolor{red}{\sqrt{h}}u  & & \text {in}\;(0,T), \\
		& \psi(1,\cdot)=0 & & \text {in}\;(0,T), \\
		&h^{\prime}+2\mu \psi_{x}(1,\cdot)=0 & & \text {in}\;(0,T),  \\	
		& h(0)=h_{0},\\
		& \psi(\cdot,0)=\psi_{0} & &\text {in}\;(0,1),
	\end{aligned}
	\right.
\end{equation}
where $h_{0}:=L_{0}^{2}$ and $\psi_{0}(x):=\varrho_{0}(xL_{0})$, for $x$ in $(0,1)$. Consequently, \eqref{s1} is locally exactly controllable to the trajectory $(\bar{\varrho},\bar{L},\bar{u})$ is equivalent that, \eqref{s2} is locally exactly controllable to the trajectory $(\overline{\psi},\overline{h},\overline{u})$ of system \eqref{s2}, where  
$(\bar{\varrho},\bar{L},\bar{u})$ and $(\overline{\psi},\overline{h},\overline{u})$ are linked in compliance with \eqref{Chang}.
\par 
Next, we will transform the local exact controllability to the trajectories of \eqref{s2} into a problem of local null controllability. To accomplish this, we will introduce the change of variable:
\begin{eqnarray*}
	z=\psi-\overline{\psi}\quad \mbox{and}\quad  k=h-\overline{h},
\end{eqnarray*}
where $(\overline{\psi},\overline{h},\overline{u})$ is a trajectory of \eqref{s2} and keep all the terms which are linear with respect to $(z,k)$. The nonlinear system satisfied by the perturbation variables is given by
\begin{equation} \label{s3}
	\left\{
	\begin{aligned}
		&\overline{h}z_{t}-z_{xx} -\frac{x}{2}(\overline{h})^{\prime}z_{x} +\mu x(\overline{\psi})_{x}z_{x}(1,t)+mk+nz =\mathcal{Q}(z,k) & & \text {in}\; Q, \\
		&z_{x}(0,t)=\hat{u}(t) & & \text {in}\;(0,T), \\
		& z(1,t)=0  & & \text {in}\;(0,T), \\
		&k^{\prime}(t)+2\mu z_{x}(1,t)=0 & & \text {in}\;(0,T),  \\	
		& k(0)=k_{0},\\
		& z(\cdot,0)=z_{0} & &\text {in}\;(0,1),
	\end{aligned}
	\right.
\end{equation}
where $k^{\prime}$ was replaced by $-2\mu z_{x}(1,t)$ in the first equation of \eqref{s3}, $\hat{u}:=\textcolor{red}{\sqrt{h}}u-\textcolor{red}{\sqrt{\overline{h}}}\overline{u}$, $k_{0}:=h_{0}-\overline{h}(0)$, $z_{0}:=\psi_{0}-\overline{\psi}(\cdot,0)$, the smooth bounded coefficients are given by
\begin{eqnarray}
	m:=(\overline{\psi})_{t}-\overline{\psi}(1-\overline{\psi}), \quad n:=-\overline{h}(1-2\overline{\psi}) \label{Coiff}
\end{eqnarray}
and all nonlinear terms have been grouped into $\mathcal{Q}(z,k)$:
\begin{eqnarray}
	\mathcal{Q}(z,k):=-kz_{t}-\mu x z_x(1,t)z_x-kz^2+(1-2\overline{\psi})kz-\overline{h}z^2. \label{nonlinear Coiff}
\end{eqnarray}
\par
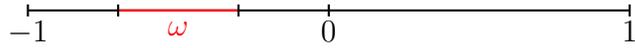
\begin{figure}[H]
	\centering
	\begin{tikzpicture}[scale=4] 
		\draw[thick] (-1,0) -- (1,0) node[below right] {};
		
		\draw[thick, red] (-0.7,0) -- (-0.3,0) node[midway, below] {$\omega$};
		\foreach \x in {-1,-0.7,-0.3,0,1} {
			\draw[thick] (\x,0.02) -- (\x,-0.02);
		}
		
		\node[below] at (-1,0) {$-1$};
		\node[below] at (0,0) {$0$};
		\node[below] at (1,0) {$1$};
	\end{tikzpicture}
	\caption{The extended domain $(-1,1)$ of $(0,1)$ and the control region $\omega$.} 
	\label{fig:ext domain}
\end{figure}

 
The strategy adopted to control system \eqref{s3} is based on the extension-restriction principle, which transforms a boundary control problem into an interior control problem. This approach involves extending the system to a larger domain where the control acts internally, outside the original domain. Once the solution is controlled in this extended domain, its restriction to the original domain yields a controlled trajectory with Neumann boundary control at $x=0$ and Dirichlet conditions at $x=1$.
Hence, we reformulate the local null controllability of \eqref{s3} as local null controllability with interior control by extending the domain $(0,1)$ into $(-1,1)$ with control acts in a region of $(-1,0)$ as shown in Figure \ref{fig:ext domain}:
\begin{equation} \label{s4}
	\left\{
	\begin{aligned}
		&\overline{h}z_{t}-z_{xx} -\frac{x}{2}(\overline{h})^{\prime}z_{x}+ \mu x (\overline{\psi})_{x}z_{x}(1,t)+mk+nz=\mathcal{Q}(z,k)+\mathds{1}_{\omega}v & & \text {in}\; \tilde{Q}, \\
		&z(-1,t)=z(1,t)=0  & & \text {in}\;(0,T), \\
		&k^{\prime}(t)+2\mu z_{x}(1,t)=0 & & \text {in}\;(0,T),  \\	
		& k(0)=k_{0},\\
		& z(\cdot,0)=z_{0} & &\text {in}\;(-1,1),
	\end{aligned}
	\right.
\end{equation}
where $\tilde{Q}:=(-1,1)\times (0,1)$ and $\omega\subset\subset (-1,0)$ is an open subset of $(-1,0)$.
\begin{remark}
	To simplify, we have retained the same notations $z_0$, $\overline{\psi}$ and $\overline{h}$ for the new extended data of system \eqref{s4}. Specifically, $z_0$ has been extended as an element of $H^{1}(-1,1)$, while $\overline{\psi}$ and $\overline{h}$ have been extended as follows:
	\begin{eqnarray}
		(\overline{\psi},\overline{h})\in \left[H^{2,1}(\tilde{Q})\cap W^{1,\infty}(0,T;H^{1}(0,T))\right]\times W^{1,\infty}(0,T), \label{SC}
	\end{eqnarray} 
	which satisfies $\overline{h}>\textcolor{red}{ L_\star^2}$, in order to ensure the regularity required for investigating the Carleman estimate of \cite{barcena2023exact}.
\end{remark}
\par 
Considering all the above, the main Theorem \ref{main result} will be an immediate consequence of the following theorem: 
\begin{theorem} \label{result}
	There exists $\varepsilon> 0$ such that for any $k_{0}\in\mathbb{R}$ and any $z_{0}\in H_{0}^{1}(-1,1)$ satisfying
	\begin{eqnarray*}
		|k_{0}|+ \|z_{0}\|_{H^{1}(-1,1)}< \varepsilon,
	\end{eqnarray*}
	there exists a control $v\in L^{2}(\omega\times(0,T))$ such that the associated state $(z, k)$ of system \eqref{s4} belonging to $H^{2,1}_{0}(\tilde{Q})\times H^{1}(0,T)$ and \textcolor{red}{$\|(z, k)\|_{H^{2,1}_{0}(\tilde{Q})\times H^{1}(0,T)} \leq C(|k_{0}|+ \|z_{0}\|_{H^{1}(-1,1)})$ for some constant $C>0$ independent of $\varepsilon$}, satisfies:
	\begin{eqnarray*}
		z(\cdot,T)=0 \;\mbox{in}\; (-1,1) \quad \mbox{and} \quad k(T)=0,
	\end{eqnarray*}
	where 
	\begin{eqnarray*}
		H_{0}^{2,1}(\tilde{Q}):=\left\{z\in H^{2,1}(\tilde{Q})\;:\; z(-1,t)=z(1,t)=0\; \mbox{in}\; (0,T)\right\}.
	\end{eqnarray*}
	\qed
\end{theorem}
\textcolor{red}{
\begin{remark} \label{reg of contr}
	By the energy estimate given in Theorem \ref{result} and since  $\overline{h}\textcolor{red}{>} L_\star^2$, we can choose $\varepsilon>0$ sufficiently small so that
	\begin{align*}
		|k(t)|< \min_{t\in [0,T]}\overline{h}(t) - L_\star^2 \quad \forall t\in [0,T].
	\end{align*}
	Consequently, we obtain
	\begin{align*}
		h(t)&=k(t)+\overline{h}(t)> L_\star^2\quad\mbox{and}\quad L(t)=\sqrt{h(t)}> L_\star \quad\forall t\in [0,T].
%		L(t)&=\sqrt{h(t)}> L_\star,  &&\forall t\in [0,T].
	\end{align*}
	We now explain why Theorem~\ref{main result} follows from Theorem~\ref{result}. 
	Indeed, by the restriction strategy and the trace theorem 
	\cite[Theorem~2.1]{lions2012non}, the control $\hat u$ is defined by
	\[
	\hat u=z_x(x,\cdot)\big|_{x=0}\in H^{1/4}(0,T),
	\]
	where $(z, k)$ is the controlled state of system \eqref{s4}. The control $u$ appearing in Theorem~\ref{main result} is then given by
	\[
	u=\frac{1}{\sqrt{h}}\Bigl(\hat u+\sqrt{\overline h}\,\overline u\Bigr).
	\]
	Since $\overline u\in H^{1/4}(0,T)$, $h,\overline h\in H^1(0,T)$ and satisfy $h(t),\overline h(t)>L_\star^2>0, \; t\in (0,T)$, 
	the functions $1/\sqrt{h}$ and $\sqrt{\overline h}$ belong to $H^1(0,T)$.
	Consequently, using the multiplication theorem in Sobolev spaces 
	\cite[Theorem~7.4]{behzadan2021multiplication}, we deduce that $u\in H^{1/4}(0,T)$. 
\end{remark} }
To achieve this, we reformulate the local null controllability property of system \eqref{s4} as the surjectivity of the mapping 
$\Lambda: \mathds{X}\longrightarrow\mathds{Y}$, where $\mathds{X}$ and $\mathds{Y}$ are appropriately chosen Hilbert spaces (see Section \ref{Sec 3}):
\begin{equation} \label{operator0}
	\begin{aligned}
		&\Lambda(z,k,v):=\left(\Lambda_{1}(z,k,v), \Lambda_{2}(z,k,v), z(\cdot,0), k(0)\right), \\
		&\Lambda_{1}(z,k,v):=\overline{h}z_{t}-z_{xx} -\frac{x}{2}(\overline{h})^{\prime}z_{x}+ \mu x (\overline{\psi})_{x}z_{x}(1,\cdot)+mk+nz-\mathcal{Q}(z,k)-\mathds{1}_{\omega}v, \\
		&\Lambda_{2}(z,k,v):=k^{\prime}-2\mu z_{x}(\cdot,t).
	\end{aligned}
\end{equation}
More precisely, Theorem \ref{result} is equivalent to 
\begin{eqnarray*}
	\exists\varepsilon>0,\;  \forall (z_0, k_{0})\in B_{H^{1}_{0}(-1,1)\times\mathbb{R}}((0,0),\varepsilon),\; \exists (z, k, v)\in\mathbb{X}\;:\; \varLambda(z, k ,v)=(0, 0,z_0, k_0).
\end{eqnarray*}
For that, we apply the Lyusternik-Graves Inverse Mapping Theorem (see Theorem \ref{Lyusternik}). \textcolor{red}{This approach to proving local null controllability has been used in several works by Fern\'andez-Cara and co-authors, notably for quasilinear heat equations and the one-phase Stefan problem \cite{barcena2023exact, fernandez2023local}. In the context of local null controllability for volume-surface reaction-diffusion equations with dynamic boundary conditions, we refer to \cite{et2025null}.}

\noindent \textbf{Structure of this paper.} The rest of the paper is structured as follows. Section \ref{Sec 2} focuses on the well-posedness of the linearized system of \eqref{s4}, as well as the adjoint system. We also prove the null controllability of the linearized system of \eqref{s4} with a source term. In Section \ref{Sec 3}, we show the local null controllability of the nonlinear system \eqref{s4}. Section \ref{Sec 4} focuses on numerical analysis, introducing recent deep learning techniques, especially Physics-Informed Neural Networks (PINNs), to simulate the behavior of the density and free boundary in the nonlinear system, as well as its linear approximation. Finally, Section \ref{Sec 5} is devoted to the conclusion.

\section{Study of the linearized system of \eqref{s4}} \label{Sec 2}
As usual, the local inversion of a nonlinear mapping relies on the surjectivity (which is sufficient) of its derivative around the point of interest. At this stage, the surjectivity of our mapping defined in \eqref{operator0} is connected to the null controllability of the inhomogeneous linearized system (around zero) of \eqref{s4}.
\subsection{Linearization}
The inhomogenuous linearized system at $(0,0)$ of \eqref{s4} is given by:
\begin{equation} \label{s5}
	\left\{
	\begin{aligned}
		&\overline{h}z_{t}-z_{xx} -\frac{x}{2}(\overline{h})^{\prime}z_{x}+ \mu x \overline{\psi}_{x}z_{x}(1,t)+mk+ nz=f_{0}+\mathds{1}_{\omega}v & & \text {in}\; \tilde{Q}, \\
		&z(-1,t)=z(1,t)=0  & & \text {in}\;(0,T), \\
		&k^{\prime}(t)+2\mu z_{x}(1,t)=g_{0} & & \text {in}\;(0,T),  \\	
		& k(0)=k_{0},\\
		& z(\cdot,0)=z_{0} & &\text {in}\;(-1,1),
	\end{aligned}
	\right.
\end{equation}
where $m$ and $n$ are defined in \eqref{Coiff}, and the functions $f_{0}$ and $g_{0}$ belong to suitable function spaces that exhibit exponential decay as $t$ approaches $T$ from the left, as will be defined more precisely later.
In the sequel, we use the following notations:
\begin{notation} We will denote 
	\begin{equation} \label{n1} 
		\begin{aligned}
			\textbf{L}_{1}(z,k)&:=\overline{h}z_{t}-z_{xx} -\frac{x}{2}(\overline{h})^{\prime}z_{x}+ \mu x \overline{\psi}_{x}z_{x}(1,\cdot)+mk+nz,\\
			\textbf{L}_{2}(z,k)&:=k^{\prime}+2\mu z_{x}(1,\cdot),\\
			\textbf{L}_{1}^{\star}(z,k)&:=-(\overline{h}z)_{t}-z_{xx}+\frac{(\overline{h})^{\prime}}{2}(xz)_{x}+ nz,\\
			\textbf{L}_{2}^{\star}(z,k)&:=-k^{\prime}+\displaystyle\int_{-1}^{1}m(x,\cdot)z(x,\cdot)\d x.
		\end{aligned}
	\end{equation}
\end{notation}
We are focused on addressing the following categories of solutions for \eqref{s5}.
\begin{definition}
	Let $(f_{0},g_{0})\in L^{2}(\tilde{Q})\times L^{2}(0,T)$, $(z_{0},k_{0})\in L^{2}(-1,1)\times\mathbb{R}$ and $v\in L^{2}(\omega\times (0,T))$. 
	\begin{enumerate}[label=(\arabic*)]
		\item A solution by transposition of \eqref{s5} is a couple of functions $(z,k)\in L^{2}(\tilde{Q})\times L^{2}(0,T)$ such that for any $(w,\theta)\in H^{2,1}_{0}(\tilde{Q})\times H^{1}(0,T)$ with $w(\cdot,T)=\theta(T)=0$, $w(-1,t)=0$ and $w(1,t)-\mu\displaystyle\int_{-1}^{1}x\overline{\psi}_{x}w(x,t)\d x-2\mu\theta(t)=0$, $t\in (0,T)$, we have 
		\begin{eqnarray*}
			&&\iint_{\tilde{Q}}z \textbf{L}_{1}^{\star}(w,\theta) \d x \d t + \int_{0}^{T}k \textbf{L}_{2}^{\star}(w,\theta) \d t -\iint_{\omega\times (0,T)}vw \d x \d t=\iint_{\tilde{Q}}f_{0}w\d x \d t \nonumber\\
			&& \quad +\int_{0}^{T}g_{0}(t)\theta(t)\d t + \overline{h}(0)\int_{-1}^{1}z_{0}(x)w(x,0)\d x + k_{0}\theta(0).  
			%\label{w1}
		\end{eqnarray*}
		\item A strong solution of \eqref{s5} is a couple of functions $(z,k)\in H^{2,1}_{0}(\tilde{Q})\times H^{1}(0,T)$ fulfilling $\eqref{s5}_{1}$ in $L^{2}(\tilde{Q})$, $\eqref{s5}_{3}$ in $L^{2}(0,T)$ and initial conditions.
	\end{enumerate}
\end{definition}
\begin{remark}
	A transposition solution is related to the solutions of adjoint system, where conditions $w(-1,t)=0$ and $w(1,t)-\mu\displaystyle\int_{-1}^{1}x\overline{\psi}_{x}w(x,t)\d x-2\mu\theta(t)=0$, $t\in (0,T)$ appear as boundary conditions of the adjoint system (see system \eqref{s6}).
\end{remark}
System \eqref{s5} has a unique strong solution in $H_{0}^{2,1}(\tilde{Q})\times H^{1}(0,T)$, for any $f_{0}\in L^{2}(\tilde{Q}), v\in L^{2}(\omega\times (0,T)), g_{0}\in L^{2}(0,T), k_{0}\in\mathbb{R}$ and any $z_{0}\in H^{1}_{0}(-1,1)$,  see Proposition 2.7 in \cite{barcena2023exact}. Moreover, there is a positive constant $C$ depending on $\overline{\psi}, \overline{h}, \mu$ and $T$ such that 
\begin{eqnarray} 
	\|z\|_{H_{0}^{2,1}(\tilde{Q})} + \|k\|_{H^{1}(0,T)}\leq C\left(\|f_{0}\|_{L^{2}(\tilde{Q})}+\|v\|_{L^{2}(\omega\times (0,T))}+ \|g_{0}\|_{L^{2}(0,T)} + |k_{0}|+ \|z_{0}\|_{H^{1}_{0}(-1,1)} \right). \label{Energy estiamte D}
\end{eqnarray}

\subsection{Adjoint system of \eqref{s5}}
To determine the adjoint system of \eqref{s5}, multiply the first equation of \eqref{s5} by $w:=w(x,t)$ and the third equation by $\theta:=\theta(t)$ and proceeding with integrations by parts, we find that the adjoint system of \eqref{s5} is given by
\begin{equation} \label{s6} 
	\left\{
	\begin{aligned}
		&-(\overline{h}w)_{t}-w_{xx}+\frac{(\overline{h})^{\prime}}{2}(xw)_{x}+ nw=f_{1} & & \text {in}\;  \tilde{Q}, \\
		&w(-1,t)=0 & & \text {in}\;(0,T), \\
		& w(1,t)=\mu\int_{-1}^{1} x \overline{\psi}_{x}(x,t) w(x,t)\d x +2\mu\theta (t)  & & \text {in}\;(0,T), \\
		&-\theta^{\prime}(t)+\int_{-1}^{1}m(x,t)w(x,t)\d x=g_{1} & & \text {in}\;(0,T),  \\	
		& \theta(T)=\theta_{T},\\
		& w(\cdot,T)=w_{T}, & &\text {in}\;(-1,1).
	\end{aligned}
	\right.
\end{equation}
Note that the data must satisfy the compatibility conditions:
\begin{eqnarray}
	&&w_{T}(-1)=0\quad \mbox{and} \quad w_{T}(1)=\mu\int_{-1}^{1} x \overline{\psi}_{x}(x,T) w_{T}(x)\d x +2\mu\theta_{T}. \label{compatibility condition A}
\end{eqnarray}
The following result assures the existence of strong solutions of \eqref{s6} and its proof is a direct application of \cite[Proposition 2.5]{barcena2023exact} based on Leray-Schauder's Fixed-Point Principle.
\begin{proposition}
	Let $w_{T}\in H^{1}(-1,1)$ and  $\theta_{T}\in\mathbb{R}$ such that the compatibility condition \eqref{compatibility condition A} is satisfied and  $(f_{1},g_{1})\in L^{2}(\tilde{Q})\times L^{2}(0,T)$, then \eqref{s6} has a unique strong solution $(w,\theta)\in H^{2,1}(\tilde{Q})\times H^{1}(0,T)$, which satisfies the following
	estimate:
	\begin{eqnarray}
		\|w\|_{H^{2,1}(\tilde{Q})} + \|\theta\|_{H^{1}(0,T)}\leq C\left(\|f_{1}\|_{L^{2}(\tilde{Q})}+ \|g_{1}\|_{L^{2}(0,T)} + |\theta_{T}|+ \|w_{T}\|_{H^{1}(-1,1)} \right), \label{Energy estiamte}
	\end{eqnarray}
	where $C$ is a positive constant depending on $\overline{h}, \overline{\psi}, \mu$ and $T$.
\end{proposition}

\subsection{Null controllability of system \eqref{s5}}
The objective of this section is to prove null controllability for system \eqref{s5} and some state estimation and control regularity that requires certain regularity of the initial data and inhomogeneous terms. We consider the following positive weight functions $\alpha$ and $\xi$, which depend on $\omega$:
$$\alpha(x, t):=\frac{e^{2\lambda \kappa}-e^{\lambda\left(\kappa+ \eta(x)\right)}}{t(T-t)}\quad \mbox{and} \quad \xi(x, t):=\frac{e^{\lambda\left(\kappa+ \eta(x)\right)}}{t(T-t)}.
$$
Here, $\lambda, \kappa >1$ and $\eta=\eta(x)$ is a function in $C^2([-1,1])$ satisfying
\begin{equation}
	\label{eta}
	\eta>0 \text { in } [-1,1],\quad  \min_{x\in [-1,1] \backslash \omega^{\prime}}\left| \eta_{x}(x)\right| >0,\quad \eta(-1)=\eta(1)=\min_{x\in [-1,1]}\eta(x),
\end{equation}
even if we divide by $||\eta||_{\infty}$, we can also assume that $||\eta||_{\infty}=1$,
where $\omega^{\prime}\subset\subset \omega$ is an open non empty subset. The existence of such a function $\eta$ satisfying \eqref{eta} is proved in \cite{fursikov1996controllability}.
\begin{notation} For all $t\in (0,T)$, we denote:
	\begin{eqnarray*}
		\hat{\alpha}(t):=\max_{x\in [-1,1]}\alpha(x,t)=\alpha(\pm 1,t)\quad \mbox{and} \quad \check{\xi}(t):=\min_{x\in [-1,1]}\xi(x, t)=\xi(\pm 1,t).
	\end{eqnarray*}
\end{notation}
\begin{proposition} \label{P2}
	There exist constants $\lambda_{1}, s_{1}\geq 1$ and $C_{1}>0$ such that, for any $\lambda\geq\lambda_{1}$, any $s\geq s_{1}(T+T^{2})$, any $w_{T}\in H^{1}(-1,1)$ and any $\theta_{T}\in\mathbb{R}$ with compatibility condition \eqref{compatibility condition A} is satisfied and  any source terms $(f_{1},g_{1})\in L^{2}(\tilde{Q})\times L^{2}(0,T)$, the strong solution to \eqref{s6} satisfies:
	\begin{eqnarray}
		&&\iint_{\tilde{Q}}e^{-2s\alpha}\left[(s\xi)^{-1}(|w_{xx}|^{2}+|w_{t}|^{2})+\lambda^{2}(s\xi)|w_{x}|^{2}+\lambda^{4}(s\xi)^{3}|w|^{2}\right]\d x\d t \nonumber \\
		&&   + \int_{0}^{T}e^{-2s\hat{\alpha}}\left[\lambda^{3}(s\check{\xi})^{3}|w(1,t)|^{2}+ \lambda(s\check{\xi})\left(|w_{x}(-1,t)|^{2}+|w_{x}(1,t)|^{2}\right) \right]\d t \nonumber \\
		&&   + \int_{0}^{T}e^{-2s\hat{\alpha}}\left[|\theta_{t}|^{2}+\lambda^{3}(s\check{\xi})^{3}|\theta|^{2}\right]\d t \nonumber \\
		&&\leq C_{1} \left(s^{3}\lambda^{4}\iint_{\omega\times (0,T)}e^{-2s\alpha}\xi^{3}|w|^{2}\d x\d t+ \iint_{\tilde{Q}}e^{-2s\alpha}|f_{1}|^{2}\d x\d t + \int_{0}^{T}e^{-2s\hat{\alpha}}|g_{1}|^{2}\d t\right). \label{Carleman 1}
	\end{eqnarray}
\end{proposition}
\begin{proof}
	System \eqref{s6} can be written in the form
	\begin{equation} \label{s7} 
		\left\{
		\begin{aligned}
			&\varphi_{t}+d(t)\varphi_{xx}=f & & \text {in}\;  \tilde{Q}, \\
			&\varphi(-1,t)=0 & & \text {in}\;(0,T), \\
			& \varphi(1,t)=\gamma(t)+\int_{-1}^{1} N(x,t) \varphi(x,t)\d x   & & \text {in}\;(0,T), \\
			&\gamma^{\prime}(t)-\int_{-1}^{1}R(x,t)\varphi(x,t)\d x=g & & \text {in}\;(0,T),  \\	
			& \gamma(T)=2\mu\theta_{T},\\
			& \varphi(\cdot,T)=\varphi_{T} & &\text {in}\;(-1,1),
		\end{aligned}
		\right.
	\end{equation}
	where $(\varphi,\gamma):=(w,2\mu\theta)$, $d:=(\overline{h})^{-1},\; N:=\mu x\overline{\psi}_{x},\; R:=2\mu m$, $f:=-d\left(f_{1}+\left(\frac{(\overline{h})^{\prime}}{2} -n\right)\varphi-\frac{x}{2}(\overline{h})^{\prime}\varphi_{x}\right)$, $g:=-2\mu g_{1}$, \textcolor{red}{$m$ and $n$ are defined in \eqref{Coiff}}. Since $(\overline{\psi},\overline{h})$ is a smooth data (see \eqref{SC}) \textcolor{red}{such that $\overline{h}> L_\star^2$}, using the Sobolev continuous embedding $H^{2,1}(\tilde{Q})\hookrightarrow C([0,T];H^{1}(-1,1))$ and $H^{1}(-1,1)$ is a Banach algebra, we obtain
	\begin{eqnarray*}
		R\in L^{\infty}(0,T;L^{2}(-1,1)),\; N\in W^{1,\infty}(0,T;L^{2}(-1,1))\;\text{and}\; d\in W^{1,\infty}(0,T).
	\end{eqnarray*}
	Note that the regularity $W^{1,\infty}(0,T)$ is sufficient for $d$. Hence, we are in a position to apply Theorem 3.1 of \cite{barcena2023exact}, so there exist constants $\lambda_{1}, s_{1}\geq 1$ and $C>0$ such that, for any $\lambda\geq\lambda_{1}$, any $s\geq s_{1}(T+T^{2})$, the strong solution to \eqref{s7} satisfies:
	\begin{eqnarray*}
		&&\iint_{\tilde{Q}}e^{-2s\alpha}\left[(s\xi)^{-1}(|\varphi_{xx}|^{2}+|\varphi_{t}|^{2})+\lambda^{2}(s\xi)|\varphi_{x}|^{2}+\lambda^{4}(s\xi)^{3}|\varphi|^{2}\right]\d x\d t\\
		&&  \quad + \int_{0}^{T}e^{-2s\hat{\alpha}}\left[\lambda^{3}(s\check{\xi})^{3}|\varphi(1,t)|^{2}+ \lambda(s\check{\xi})\left(|\varphi_{x}(-1,t)|^{2}+|\varphi_{x}(1,t)|^{2}\right) \right]\d t \\
		&&  \quad + \int_{0}^{T}e^{-2s\hat{\alpha}}\left[|\gamma_{t}|^{2}+\lambda^{3}(s\check{\xi})^{3}|\gamma|^{2}\right]\d t \\
		&&\quad\leq C \left(s^{3}\lambda^{4}\iint_{\omega\times (0,T)}e^{-2s\alpha}\xi^{3}|\varphi|^{2}\d x\d t+ \iint_{\tilde{Q}}e^{-2s\alpha}|f|^{2}\d x\d t+ \int_{0}^{T}e^{-2s\hat{\alpha}}|g|^{2}\d t\right).
	\end{eqnarray*}
	Absorbing the lower-order terms of $f$, we easily obtain the required inequality. 
\end{proof}
We will deduce a Carleman estimate with functions blowing up only at $t=T$, which allows us to
prove an observability inequality for system \eqref{s6}.
Define the new weight functions:
\begin{eqnarray*}
	\beta(x, t):=\frac{e^{2\lambda \kappa}-e^{\lambda\left(\kappa+ \eta(x)\right)}}{\ell(t)}\quad \mbox{and} \quad \zeta(x, t):=\frac{e^{\lambda\left(\kappa+ \eta(x)\right)}}{\ell(t)},\quad (x,t)\in [-1,1]\times (0,T),
\end{eqnarray*}
where the function $\ell$ is given by
\begin{eqnarray*}
	\ell(t):=\begin{cases}  \frac{T^{2}}{4}
		\quad & \mbox{if}\; t\in [0,T/2),\\
		t(T-t)   \quad & \mbox{if}\; t\in [T/2,T].
	\end{cases}
\end{eqnarray*}
Note that $\ell\in C^{1}([0,T])$.
An estimate with such weights is given in the following result. In the proof, we will use Proposition \ref{P2} and energy estimates \eqref{Energy estiamte}.
\begin{proposition}
	There exist constants $\lambda_1, s_{1}\geq 1$ such that for any $s\geq s_1$, any $\lambda\geq\lambda_1$, there exists a constant $C:=C(s,\lambda,T)$, which satisfies the following property: for any strong solution $(w,\theta)$ of \eqref{s6}, we have the following estimate:
	\begin{eqnarray}
		&&\iint_{\tilde{Q}}e^{-2s\beta}\left[\ell(|w_{xx}|^{2}+|w_{t}|^{2})+\ell^{-1}|w_{x}|^{2}+\ell^{-3}|w|^{2}\right]\d x\d t \nonumber\\
		&&   + \int_{0}^{T}e^{-2s\hat{\beta}}\left[\ell^{-3}|w(1,t)|^{2}+ \ell^{-1}\left(|w_{x}(-1,t)|^{2}+|w_{x}(1,t)|^{2}\right) \right]\d t \nonumber \\
		&&   + \int_{0}^{T}e^{-2s\hat{\beta}}\left[|\theta_{t}|^{2}+\ell^{-3}|\theta|^{2}\right]\d t + \int_{-1}^{1}|w(x,0)|^{2}\d x+ |\theta(0)|^{2} \nonumber \\
		&&\leq C \left(\iint_{\omega\times (0,T)}e^{-2s\beta}\ell^{-3}|w|^{2}\d x\d t+ \iint_{\tilde{Q}}e^{-2s\beta}|f_{1}|^{2}\d x\d t+ \int_{0}^{T}e^{-2s\hat{\beta}}|g_{1}|^{2}\d t\right). \label{Carleman 2}
	\end{eqnarray}
	with $s_{1}$ and $\lambda_{1}$ as in Proposition \ref{P2} and $\hat{\beta}=\displaystyle\max_{x\in [-1,1]}\beta(x,\cdot)$.
\end{proposition}
\begin{proof}
	To do so, we split the left-hand side of \eqref{Carleman 2} into two parts, corresponding to the restrictions of $(w,\theta)$ on $[-1,1]\times\left(0,T/2\right)$ and $[-1,1]\times\left(T/2,T\right)$.\\
	Firstly, on $[-1,1]\times\left(T/2,T\right)$. Using Carleman's estimate \eqref{Carleman 1}, $\alpha=\beta$ and $\xi=e^{\lambda(\kappa+\eta(x))}\ell^{-1}$ in $[-1,1]\times \left(T/2,T\right)$, we obtain 
	\begin{eqnarray}
		&&\int_{-1}^{1}\int_{T/2}^{T}e^{-2s\beta}\left[\ell(|w_{xx}|^{2}+|w_{t}|^{2})+\ell^{-1}|w_{x}|^{2}+\ell^{-3}|w|^{2}\right]\d x\d t \nonumber\\
		&&   + \int_{T/2}^{T}e^{-2s\hat{\beta}}\left[\ell^{-3}|w(1,t)|^{2}+ \ell^{-1}\left(|w_{x}(-1,t)|^{2}+|w_{x}(1,t)|^{2}\right) \right]\d t \nonumber\\
		&&   + \int_{T/2}^{T}e^{-2s\hat{\beta}}\left[|\theta_{t}|^{2}+\ell^{-3}|\theta|^{2}\right]\d t \nonumber \\
		&&\leq C \left(\iint_{\omega\times (0,T)}e^{-2s\alpha}\xi^{3}|w|^{2}\d x\d t+ \iint_{\tilde{Q}}e^{-2s\alpha}|f_{1}|^{2}\d x\d t + \int_{0}^{T}e^{-2s\hat{\alpha}}|g_{1}|^{2}\d t\right). \label{c2}
	\end{eqnarray}
	Using the definition of $\ell$, we have $\ell(t)\geq t(T-t)$ on $[0,T]$, then
	\begin{eqnarray}
		\alpha\geq \beta \quad \mbox{in}\;\; [-1,1]\times (0,T) \quad \mbox{and} \quad \hat{\alpha}\geq \hat{\beta} \quad \mbox{in}\;\; (0,T). \label{c3}
	\end{eqnarray}
	On the other hand, one has
	\begin{eqnarray}
		e^{-2s\alpha}\xi^{3}\leq C e^{-2s\beta}\ell^{-3}\quad \mbox{in}\;\; [-1,1]\times (0,T), \label{c4}
	\end{eqnarray}
	where $C>0$ only depending on $s$, $\lambda$ and $T$. 
	%Indeed, in $[-1,1]\times (T/2,T)$, we have $e^{-2s\alpha}\xi^{3}= e^{3\lambda(\kappa+\eta(x))}e^{-2s\beta}\ell^{-3}$ and in $[-1,1]\times (0,T/2)$,  $e^{-2s(\alpha-\beta)}\left(\xi\ell\right)^{3}$ is bounded with respect to $x$ and this bound admits $0$ as a limit when $t\searrow 0$, hence $e^{-2s(\alpha-\beta)}\left(\xi\ell\right)^{3}$ is bounded on $[-1,1]\times (0,T/2)$.\\
	Considering \eqref{c2}-\eqref{c4}, we obtain 
	\begin{eqnarray*}
		&&\int_{-1}^{1}\int_{T/2}^{T}e^{-2s\beta}\left[\ell(|w_{xx}|^{2}+|w_{t}|^{2})+\ell^{-1}|w_{x}|^{2}+\ell^{-3}|w|^{2}\right]\d x\d t\\
		&&   + \int_{T/2}^{T}e^{-2s\hat{\beta}}\left[\ell^{-3}|w(1,t)|^{2}+ \ell^{-1}\left(|w_{x}(-1,t)|^{2}+|w_{x}(1,t)|^{2}\right) \right]\d t \\
		&&   + \int_{T/2}^{T}e^{-2s\hat{\beta}}\left[|\theta_{t}|^{2}+\ell^{-3}|\theta|^{2}\right]\d t \\
		&&\leq C \left(\iint_{\omega\times (0,T)}e^{-2s\beta}\ell^{-3}|w|^{2}\d x\d t+ \iint_{\tilde{Q}}e^{-2s\beta}|f_{1}|^{2}\d x\d t+ \int_{0}^{T}e^{-2s\hat{\beta}}|g_{1}|^{2}\d t\right).
	\end{eqnarray*}
	Secondly, on $[-1,1]\times\left(0,T/2\right)$.
	Let us define the function $\vartheta\in C^{1}([0,T])$ such that 
	\begin{eqnarray*}
		\vartheta(t)= 1\quad\mbox{in}\quad \left[0,T/2\right] \quad \mbox{and}\quad \vartheta(t)= 0\quad\mbox{in}\quad \left[3T/4,T\right].
	\end{eqnarray*}
	Put  $(w^{\vartheta},\theta^{\vartheta}):=(\vartheta w,\vartheta\theta)$, then $(w^{\vartheta},\theta^{\vartheta})$ is the strong solution of the system:
	\begin{equation*} 
		\left\{
		\begin{aligned}
			&-(\overline{h}w^{\vartheta})_{t}-w^{\vartheta}_{xx}+\frac{(\overline{h})^{\prime}}{2}(xw^{\vartheta})_{x}+ nw^{\vartheta}=f^{\vartheta} & & \text {in}\;  \tilde{Q}, \\
			&w^{\vartheta}(-1,t)=0 & & \text {in}\;(0,T), \\
			& w^{\vartheta}(1,t)=\mu\int_{-1}^{1} x \overline{\psi}_{x}(x,t) w^{\vartheta}(x,t)\d x +2\mu\theta^{\vartheta}(t)  & & \text {in}\;(0,T), \\
			&-\theta^{\vartheta}_{t}(t)+\int_{-1}^{1}m(x,t)w^{\vartheta}(x,t)\d x=g^{\vartheta} & & \text {in}\;(0,T),  \\	
			& \theta^{\vartheta}(T)=0,\\
			& w^{\vartheta}(\cdot,T)=0 & &\text {in}\;(-1,1),
		\end{aligned}
		\right.
	\end{equation*}
	where $(f^{\vartheta},g^{\vartheta}):=\vartheta (f_{1},g_{1})-\vartheta_{t}(\overline{h}w,\theta)\in L^{2}(\tilde{Q})\times L^{2}(0,T)$.
	%$f^{\vartheta}:=\vartheta f_{1}-\vartheta_{t}\overline{h}w$ and $g^{\vartheta}:=\vartheta g_{1}-\vartheta_{t}\theta$.
	From energy estimate \eqref{Energy estiamte}, one has  
	\begin{eqnarray*}
		\|w^{\vartheta}\|_{H^{2,1}(\tilde{Q})} + \|\theta^{\vartheta}\|_{H^{1}(0,T)}\leq C\left(\|f^{\vartheta}\|_{L^{2}(\tilde{Q})}+ \|g^{\vartheta}\|_{L^{2}(0,T)} \right). 
	\end{eqnarray*}
	Consequently, the definition of $\vartheta$ implies that  
	\begin{eqnarray}
		&&\|w\|_{H^{2,1}((-1,1)\times (0,T/2))} + \|\theta\|_{H^{1}(0,T/2)}\leq C\left(\|f_{1}\|_{L^{2}((-1,1)\times (0,3T/4))} \right. \nonumber \\
		&&\quad \left. + \|g_{1}\|_{L^{2}(0,3T/4)} +\|w\|_{L^{2}((-1,1)\times (T/2,3T/4))}+ \|\theta\|_{L^{2}(T/2,3T/4)} \right). \label{c5}
	\end{eqnarray}
	Using the boundedness from above and from below of the weight functions $\beta$, $\hat{\beta}$ and $\ell$ in $[-1,1]\times\left(0,3T/4\right)$ and \eqref{c5}, we obtain
	\begin{eqnarray*}
		&&\int_{-1}^{1}\int_{0}^{T/2}e^{-2s\beta}\left[\ell(|w_{xx}|^{2}+|w_{t}|^{2})+\ell^{-1}|w_{x}|^{2}+\ell^{-3}|w|^{2}\right]\d x\d t\\
		&&   + \int_{0}^{T/2}e^{-2s\hat{\beta}}\left[\ell^{-3}|w(1,t)|^{2}+ \ell^{-1}\left(|w_{x}(-1,t)|^{2}+|w_{x}(1,t)|^{2}\right) \right]\d t \\
		&&   + \int_{0}^{T/2}e^{-2s\hat{\beta}}\left[|\theta_{t}|^{2}+\ell^{-3}|\theta|^{2}\right]\d t +\int_{-1}^{1}|w(x,0)|^{2}\d x+ |\theta(0)|^{2}\\
		&& \leq C\left( \int_{-1}^{1}\int_{0}^{3T/4}e^{-2s\beta}|f_{1}|^{2}\d x\d t + \int_{0}^{3T/4}e^{-2s\hat{\beta}}|g_{1}|^{2}\d t\right. \nonumber \\
		&& + \left. \int_{-1}^{1}\int_{T/2}^{3T/4} e^{-2s\beta}\ell^{-3}|w|^{2}\d x\d t + \int_{T/2}^{3T/4} e^{-2s\hat{\beta}}\ell^{-3}|\theta|^{2}\d t\right),\nonumber
	\end{eqnarray*}
	which, combined with \eqref{c5}, yields \eqref{Carleman 2}.
\end{proof}
From now on, we fix $s=s_{1}$ and $\lambda=\lambda_{1}$ and introduce the following weights, which we will need in the sequel. 

\begin{remark} \label{Remark 1}
	Taking $\kappa> 1$ large enough for instance $\kappa> \frac{\log(2e^{\lambda}-1)}{\lambda}$, we have that
	\begin{eqnarray*}
		\hat{\beta}(t)< 2 \check{\beta}(t) \quad \forall (x,t)\in [-1,1]\times (0,T),
	\end{eqnarray*}
	where
	\begin{eqnarray*}
		\hat{\beta}(t):=\max_{x\in [-1,1]}\beta (x,t),\quad \check{\beta}(t):=\min_{x\in [-1,1]}\beta (x,t).
	\end{eqnarray*}
\end{remark} 
\begin{notation} Let us introduce the notations:
	\begin{eqnarray*}
		\rho(t):=e^{s\hat{\beta}}\ell^{3/2}, \quad \rho_{0}(t):=(e^{s\hat{\beta}})^{1/2}, \quad \rho_{1}(t):=\rho_{0}(t)\ell^{3/2}, \quad \rho_{2}(t):=(e^{s\hat{\beta}})^{1/3}.
	\end{eqnarray*}
\end{notation}

\begin{proposition} \label{P1} Let $(z_{0},k_{0})\in H^{1}_{0}(-1,1)\times\mathbb{R}$ and $(f_{0},g_{0})$ are source terms such that $\rho(f_{0},g_{0})\in L^{2}(\tilde{Q})\times L^{2}(0,T)$.
	Then, there exists a control $v^{*}$ such that the strong solution $(z^{*},k^{*})$ of \eqref{s5} corresponding to $v^{*}$, $(f_{0},g_{0})$ and $(z_{0},k_{0})$, satisfies:
	\begin{eqnarray}
		&&\|\rho_{0}z^{*}\|_{L^{2}(\tilde{Q})}^{2} + \|\rho_{0}k^{*}\|_{L^{2}(0,T)}^{2}+ \|\rho_{1}v^{*}\|^{2}_{L^{2}(\omega\times (0,T))} \nonumber \\
		&&\quad\quad \leq C\left (\|\rho f_{0}\|_{L^{2}(\tilde{Q})}^{2}+ \|\rho g_{0}\|_{L^{2}(0,T)}^{2}+ \|z_{0}\|_{L^{2}(0,T)}^{2} + |k_{0}|^{2}\right). \label{c21}
	\end{eqnarray}
	In particular \eqref{s5} is null controllable. Moreover, $\rho_{2}(z^{*},k^{*})\in H^{2,1}(\tilde{Q})\times H^{1}(0,T)$  
	\begin{eqnarray}
		\|\rho_{2}z^{*}\|^{2}_{H^{2,1}(\tilde{Q})} + \|\rho_{2}k^{*}\|^{2}_{H^{1}(0,T)}\leq C\left (\|\rho f_{0}\|_{L^{2}(\tilde{Q})}^{2}+ \|\rho g_{0}\|_{L^{2}(0,T)}^{2}+ \|z_{0}\|_{L^{2}(0,T)}^{2} + |k_{0}|^{2}\right).  \label{c22}
	\end{eqnarray}
\end{proposition}
\begin{proof} 
	The proof of this result is inspired by the method of Fursikov and Imanuvilov  \cite{fursikov1996controllability}. Let us consider the following space:
	\begin{eqnarray*}
		&&\mathbf{P}:=\left\{(w,\theta)\in C^{2}([-1,1]\times [0,T])\times C^{1}([0,T])\;:\; w(-1,t)=0 \;\mbox{and}\; \right.\\
		&& \left.  \quad\quad\quad\quad w(1,t)=\mu\int_{-1}^{1} x \overline{\psi}_{x}(x,t) w(x,t)\d x +2\mu \theta(t) \quad \mbox{in}\; (0,T)
		\right\}.
	\end{eqnarray*}
	We define the bilinear form $\mathbf{B}: \mathbf{P}\times \mathbf{P}\longrightarrow\mathbb{R}$ by
	\begin{align*}
		\mathbf{B} ((w,\theta),(\hat{w},\hat{\theta}))&:=\iint_{\tilde{Q}}\rho_{0}^{-2}\textbf{L}_{1}^{\star}(w,\theta) \textbf{L}^{\star}_{1}(\hat{w},\hat{\theta})\d x \d t + \int_{0}^{T}\rho_{0}^{-2}\textbf{L}_{2}^{\star}(w,\theta) \textbf{L}^{\star}_{2}(\hat{w},\hat{\theta})\d t \\
		& + \iint_{\omega\times (0,T)}\rho_{1}^{-2} w\hat{w} \d x \d t,
	\end{align*}
	where 
	$\textbf{L}_{1}^{\star}$ and $\textbf{L}_{2}^{\star}$ are defined in \eqref{n1}. We also define the linear form $\mathbf{F}: \mathbf{P}\longrightarrow\mathbb{R}$ by
	\begin{align*}
		\mathbf{F}(w,\theta):=\iint_{\tilde{Q}}f_{0}w\d x \d t + \int_{0}^{T}g_{0}\theta\d t+ \overline{h}(0)\int_{-1}^{1}z_{0}(x)w(x,0)\d x + k_{0}\theta(0).
	\end{align*}
	We claim that, $\mathbf{B}$ is a scalar
	product in $\mathbf{P}$ and $\mathbf{F}$ is  continuous for the norm $\|\cdot\|_{\mathbf{B}}$ associated with the scalar product $\mathbf{B}$. Indeed, due to Carleman's estimate \eqref{Carleman 2} and Remark \ref{Remark 1}, there exists a constant $C:=C(s,\lambda,T)>0$ such that for all $(w,\theta)\in \mathbf{P}$, one has
	\begin{eqnarray}
		\iint_{\tilde{Q}}\rho^{-2}|w|^{2}\d x\d  t
		+ \int_{0}^{T}\rho^{-2}|\theta|^{2}\d t + \int_{-1}^{1}|w(x,0)|^{2}\d x+ |\theta(0)|^{2}  \leq C\,  \mathbf{P} ((w,\theta),(w,\theta)). \label{c9}
	\end{eqnarray}
	In particular, $\mathbf{B}$ is a scalar
	product in $\mathbf{P}$. To ensure the continuity of $\mathbf{F}$, using the Cauchy-Schwarz inequality and  \eqref{c9}, we obtain
	\begin{eqnarray}
		|\mathbf{F}(w,\theta)|
		\leq C \left(\|\rho f_{0} \|_{L^{2}(\tilde{Q})} + \|\rho g_{0} \|_{L^{2}(0,T)}+\|z_{0}\|_{L^{2}(-1,1)}+ |k_{0}|\right) \|(w,\theta)\|_{\mathbf{B}}. \label{c10}
	\end{eqnarray}
	In the sequel, we will denote by $\overline{\mathbf{P}}$ the completion of $\mathbf{P}$ for the norm $\|\cdot\|_{\mathbf{B}}$ and we will still denote $\mathbf{B}$ and $\mathbf{F}$ the corresponding continuous extensions.\\
	We deduce from Riesz Representation theorem that there exists a unique
	$(w^{*},\theta^{*})\in \overline{\mathbf{P}}$
	such that
	\begin{eqnarray}
		\mathbf{B}(w^{*},\theta^{*}), (w,\theta))=\mathbf{F}(w,\theta)\quad \forall (w,\theta)\in \overline{\mathbf{P}}. \label{c11}
	\end{eqnarray}
	Using \eqref{c10} and \eqref{c11}, we obtain
	\begin{eqnarray}
		\|(w^{*},\theta^{*})\|_{\mathbf{B}}
		\leq   C \left(\|\rho f_{0} \|_{L^{2}(\tilde{Q})} + \|\rho g_{0} \|_{L^{2}(0,T)}+\|z_{0}\|_{L^{2}(-1,1)}+ |k_{0}|\right). \label{c15}
	\end{eqnarray}
	Let us introduce $(z^{*},k^{*}, v^{*})$ with 
	\begin{eqnarray}
		(z^{*}, k^{*}):=\rho_{0}^{-2}(\mathbf{L}_{1}^{\star}(w^{*},\theta^{*}), \mathbf{L}_{2}^{\star}(w^{*},\theta^{*})),\quad v^{*}:=-\rho_{1}^{-2}w^{*}\mathds{1}_{\omega}. \label{c16}
	\end{eqnarray}
	According to \eqref{c16} and the definition of $\mathbf{B}$, we obtain 
	\begin{eqnarray}
		&&\iint_{\tilde{Q}}\rho_{0}^{2}|z^{*}|^{2}\d x \d t + \int_{0}^{T}\rho_{0}^{2}|k^{*}|^{2}\d t + \iint_{\omega\times (0,T)}\rho_{1}^{2}|v^{*}|^{2} \d x \d t=\|(w^{*},\theta^{*})\|^{2}_{\mathbf{B}}. \nonumber
	\end{eqnarray} 
	Using this equation and estimate \eqref{c15}, we deduce estimate \eqref{c21}. As a consequence $(z^{*},k^{*})\in L^{2}(\tilde{Q})\times L^{2}(0,T)$,  $v^{*}\in L^{2}(\omega\times (0,T))$ and from \eqref{c11}, $(z^{*},k^{*})$ is the unique solution by transposition of \eqref{s5} with the control $v^{*}$. Moreover, due to the uniqueness of solutions by transposition of \eqref{s5} and the regularity of initial data and term sources, we obtain $(z^{*},k^{*})$ is the unique strong solution of \eqref{s5}.\\
	Now, it is a matter of proving estimate \eqref{c22}, to do this, let us set 
	$(z,k):=\rho_{2}(z^{*},k^{*})$. Then $(z,k)$ is the solution of system: 
	\begin{equation} \label{s8}
		\left\{
		\begin{aligned}
			&\mathbf{L}_{1}(z,k) =\overline{h}\rho_{2,t}z^{*}+\rho_{2}f_{0}+\mathds{1}_{\omega}\rho_{2}v^{*}:=f^{*} & & \text {in}\; \tilde{Q}, \\
			&z(-1,t)=z(1,t)=0  & & \text {in}\;(0,T), \\
			&\mathbf{L}_{2}(z,k)=\rho_{2,t}k^{*}+\rho_{2}g_{0}:=g^{*} & & \text {in}\;(0,T),  \\	
			& k(0)=\rho_{2}(0)k_{0},\\
			& z(\cdot,0)=\rho_{2}(0)z_{0} & &\text {in}\;(-1,1),
		\end{aligned}
		\right.
	\end{equation}
	Using the elementary estimates:
	\begin{eqnarray}
		|\rho_{2,t}|\leq C\rho_{0},\quad \rho_{2}\leq C \rho_{1},\quad \rho_{1}\leq C \rho_{0} \label{elementary estimates}
	\end{eqnarray}
	and estimate \eqref{c21}, we obtain $(f^{*},g^{*})\in L^{2}(\tilde{Q})\times L^{2}(0,T)$. Given that $(z_{0},k_{0})\in H^{1}_{0}(-1,1)\times \mathbb{R}$, then $(z,k)\in H_{0}^{2,1}(\tilde{Q})\times H^{1}(0,T)$ is the strong solution of \eqref{s8}.
	Finally, from \eqref{Energy estiamte D}, \eqref{elementary estimates} and \eqref{c21}, we obtain \eqref{c22}.
\end{proof}

\begin{remark}
	The proof of Proposition \ref{P1} allows for identifying a null control of system \eqref{s5} by solving a variational problem, providing a numerical approach based on finite elements \cite{fernandez2013strong, fernandez2017theoretical,fernandez2023local,huaman2023local}. However, approximating the space \textbf{$P$} is highly challenging and computationally expensive. To overcome this difficulty, we chose to solve the controllability problem using PINNs (see Section \ref{Sec 4}), which enables direct handling of the nonlinear problem, unlike classical methods such as finite differences or finite elements.
\end{remark}

\section{Study of the nonlinear system \eqref{s4}} \label{Sec 3}
We will establish the local null controllability of system \eqref{s4} using \textbf{Lyusternik-Graves'} Theorem:
\begin{theorem}[\textbf{Lyusternik-Graves' Theorem}]
	\label{Lyusternik} Let $\mathds{X}$ and $\mathds{Y}$ be Banach spaces, and let $\Lambda: B_\mathds{X}(0,r)\subset \mathds{X}\rightarrow \mathds{Y}$ be a continuously differentiable mapping ($B_\mathds{X}(0,r)$ denote the ball in the space $\mathds{X}$ centered at $0$ with radius $r$). Let us assume that the derivative $\Lambda^{\prime}(0): \mathds{X}\rightarrow \mathds{Y}$ is onto and let
	us set $\xi_{0}=\Lambda(0)$. Then there exist $\varepsilon> 0$, a mapping
	$W: B_\mathds{Y}(\xi_{0},\varepsilon)\subset \mathds{Y}\rightarrow \mathds{X}$ and a constant $C> 0$ satisfying:
	\begin{itemize}
		\item $W(z)\in B_\mathds{X}(0,r)$ and $\varLambda\circ W(z)=z\quad \forall z\in B_\mathds{Y}(\xi_{0},\varepsilon)$,
		\item $\|W(z)\|_{\mathds{X}}\leq C\|z-\xi_{0}\|_{\mathds{Y}} \quad \forall z\in B_\mathds{Y}(\xi_{0},\varepsilon)$.
	\end{itemize}
\end{theorem}
We then introduce the following appropriate spaces and mapping so that the conditions of Theorem \ref{Lyusternik} are satisfied:
\begin{eqnarray*}
	&&\mathds{X}:=\left\{(z, k, v)\;: \quad\rho_{0}(z,k)\in L^{2}(\tilde{Q})\times L^{2}(0,T),\; \rho_{1}v\in L^{2}(\omega\times (0,T)),  \right.\\
	&& \left. \quad\quad\quad \rho(\mathbf{L}_{1}(z,k)-\mathds{1}_{\omega}v,\mathbf{L}_{2}(z,k))\in L^{2}(\tilde{Q})\times L^{2}(0,T),\; \rho_{2}z\in H^{2,1}(\tilde{Q}),\; \rho_{2}k\in H^{1}(0,T) \right\},\\
	&& 	\mathds{G}:=\{(f,g)\;:\; \rho (f,g)\in L^{2}(\tilde{Q})\times L^{2}(0,T)\},\nonumber\\
	&& \mathds{Y}:=\mathds{G}\times H^{1}_{0}(-1,1)\times\mathbb{R},
\end{eqnarray*}
where $\mathbf{L}_{1}$ and $\mathbf{L}_{2}$ are defined in \eqref{n1}. These spaces are naturally equipped with the following norms:
\begin{align*}
	\|(z,k,v)\|_{\mathds{X}}&:=\left(\|\rho_{0}z\|^{2}_{L^{2}(\tilde{Q})}+ \|\rho_{0}k\|^{2}_{L^{2}(0,T)}+ \|\rho_{1}v\|^{2}_{L^{2}(\omega\times (0,T)} +\|\rho(\mathbf{L}_{1}(z,k)-\mathds{1}_{\omega}v)\|^{2}_{L^{2}(\tilde{Q})} \right.\\
	& \left.  + \|\rho\mathbf{L}_{2}(z,k)\|^{2}_{L^{2}(0,T)} + \|\rho_{2}z\|^{2}_{L^{2}(\tilde{Q})}+\|\rho_{2}z_{t}\|^{2}_{L^{2}(\tilde{Q})}+ \|\rho_{2}z_{x}\|^{2}_{L^{2}(\tilde{Q})}+\|\rho_{2}z_{xx}\|^{2}_{L^{2}(\tilde{Q})} \right. \\
	& \left.  +\|\rho_{2}k\|^{2}_{L^{2}(0,T)}+\|\rho_{2}k_{t}\|^{2}_{L^{2}(0,T)} \right)^{1/2},\\
	\|(f,g)\|_{\mathbb{G}}&:=\left(\|\rho f\|^{2}_{L^{2}(\tilde{Q})}+\|\rho g\|^{2}_{L^{2}(0,T)}\right)^{1/2}.
\end{align*}
We can easily prove that these are Hilbert spaces. Let us recall the mapping $\Lambda$ defined on $\mathds{X}$ by
\begin{eqnarray*}
	\Lambda(z,k,v):=\left(\Lambda_{1}(z,k,v), \Lambda_{2}(z,k,v), z(\cdot,0), k(0)\right), %\label{operator}
\end{eqnarray*}
where
\begin{align*}
	\Lambda_{1}(z,k,v)&:= \mathbf{L}_{1}(z,k)-\mathds{1}_{\omega}v-\mathcal{Q}(z,k)\\
	&= \mathbf{L}_{1}(z,k)-\mathds{1}_{\omega}v+kz_{t}+ \mu x z_{x}(1,\cdot)z_{x}+(1-2\overline{\psi})kz +kz^{2} +\overline{h}z^{2},\\
	\Lambda_{2}(z,k,v)&:= k^{\prime}+2\mu z_{x}(1,\cdot)=\mathbf{L}_{2}(z,k).
\end{align*}
The following lemma shows that the mapping $\Lambda$ is well defined.
\begin{lemma} \label{Lemma3}
	The mapping $\Lambda : \mathds{X}\rightarrow\mathds{Y}$ is well defined and continuously differentiable. Moreover, $\Lambda^{\prime}(0,0,0)$ is given by 
	\begin{eqnarray*}
		\Lambda^{\prime}(0,0,0)(z,k,v)=(\mathbf{L}_{1}(z,k)-\mathds{1}_{\omega}v, \mathbf{L}_{2}(z,k), z(\cdot,0), k(0)) \quad \forall (z,k,v)\in\mathds{X}.
	\end{eqnarray*}
\end{lemma}
\begin{proof}
	Isolating the linear and multi-linear parts of $\Lambda$, we can write 
	\begin{eqnarray*}
		\Lambda(z,k,v)=\mathds{L}(z,k,v)+ \mathds{B}_{1} ((z,k,v),(z,k,v)) + \mathds{B}_{2}((z,k,v),(z,k,v),(z,k,v)),
	\end{eqnarray*} 
	for all $(z,k,v)\in \mathds{X}$ and for any $(z,k,v),(\tilde{z},\tilde{k},\tilde{v}),(\hat{z},\hat{k},\hat{v})\in \mathds{X}$,
	\begin{eqnarray*}
		&&\mathds{L}(z,k,v):= (\mathbf{L}_{1}(z,k)-\mathds{1}_{\omega}v, \mathbf{L}_{2}(z,k), z(\cdot,0), k(0)),\\
		&& \mathds{B}_{1}((z,k,v),(\tilde{z},\tilde{k},\tilde{v})):=(k\tilde{z}_{t}+\mu xz_{x}(1,\cdot)\tilde{z}_{x}+(1-2\overline{\psi})k\tilde{z} +\overline{h}z\tilde{z}, 0, 0, 0),\\
		&& \mathds{B}_{2}((z,k,v),(\tilde{z},\tilde{k},\tilde{v}),(\hat{z},\hat{k},\hat{v})):=(k\tilde{z}\hat{z}, 0, 0, 0).
	\end{eqnarray*}
	It suffices to show that $\mathds{L}, \mathds{B}_{1}$ and $\mathds{B}_{2}$ are bounded.
	Based on the definition of $\mathds{X}, \mathds{Y}$, the mapping $u\in H^{2,1}(\tilde{Q})\longrightarrow u(\cdot,0)\in H^{1}(-1,1)$ is bounded and the continuous embedding $H^{1}(0,T)\hookrightarrow L^{\infty}(0,T)$, we obtain that the linear mapping $\mathds{L}$
	is bounded form $\mathds{X}$ to $\mathds{Y}$. Concerning the bilinear mapping $\mathds{B}_{1}$, one has 
	\begin{eqnarray}
		&&\|\mathds{B}_{1}((z,k,v),(\tilde{z},\tilde{k},\tilde{v}))\|^{2}_{\mathds{Y}} \leq C\left(\iint_{\tilde{Q}}\rho^{2}(t)
		|k(t)|^{2}(|\tilde{z}_{t}(x,t)|^{2}+|\tilde{z}(x,t)|^{2})\d x\d t\right. \nonumber\\
		&& + \left. \iint_{\tilde{Q}}\rho^{2}(t)
		|z_{x}(1,t)|^{2}|\tilde{z}_{x}(x,t)|^{2}\d x\d t  +\iint_{\tilde{Q}}\rho^{2}(t)
		|z(x,t)|^{2}|\tilde{z}(x,t)|^{2}\d x\d t\right). \label{c6}
	\end{eqnarray}
	We apply the fact that $\rho_{2}k\in H^{1}(0,T)$ and the continuous embedding $H^{1}(0,T)\hookrightarrow L^{\infty}(0,T)$ for the first integral on the right-hand side of \eqref{c6}, the Sobolev continuous embedding $ H^{2,1}(\tilde{Q})\hookrightarrow C([0,T];H^{1}(-1,1))$, the continuity of the normal derivative:
	$|z_{x}(1,t)|\leq C\|z(\cdot,t)\|_{H^{2}(-1,1)}$ for the second one and the fact that $H^{1}(-1,1)$ is a Banach algebra for third one, we obtain 
	\begin{eqnarray*}
		\|\mathds{B
		}_{1}((z,k,v),(\tilde{z},\tilde{k},\tilde{v}))\|^{2}_{\mathds{Y}} 
		\leq C\|(z,k,v)\|^{2}_{\mathds{X}}\|(\tilde{z},\tilde{k},\tilde{v})\|^{2}_{\mathds{X}},
	\end{eqnarray*}
	which ensures that $\mathds{B}_{1}$ is bounded from $\mathds{X}$ to $\mathds{Y}$. For the multilinear mapping $\mathds{B}_{2}$, using the Sobolev continuous embeddings $H^{1}(0,T)\hookrightarrow L^{\infty}(0,T)$, $H^{1}(-1,1)$ is a Banach algebra and the Sobolev continuous embedding $ H^{2,1}(\tilde{Q})\hookrightarrow C([0,T];H^{1}(-1,1))$, we obtain
	\begin{eqnarray*}
		\|\mathds{B}_{2}((z,k,v),(\tilde{z},\tilde{k},\tilde{v}),(\hat{z},\hat{k},\hat{v}))\|^{2}_{\mathds{Y}} 
		\leq C\|(z,k,v)\|^{2}_{\mathds{X}}\|(\tilde{z},\tilde{k},\tilde{v})\|^{2}_{\mathds{X}}\|(\hat{z},\hat{k},\hat{v})\|^{2}_{\mathds{X}}.
	\end{eqnarray*}
	Hence, $\mathds{B}_{2}$ is bounded from $\mathds{X}$ to $\mathds{Y}$.
	Consequently, $\mathds{L}, \mathds{B}_{1}$ and  $\mathds{B}_{2}$ are
	continuously differentiable. In particular, $\Lambda^{\prime}(0,0,0)$ is given by 
	\begin{eqnarray*}
		\Lambda^{\prime}(0,0,0)(z,k,v)=\mathds{L}(z,k,k)\quad \forall (z,k,v)\in\mathds{X}. %\label{Dif}
	\end{eqnarray*}
\end{proof}
Due to this lemma and null controllability result of inhomogeneous linearized system \eqref{s5} obtained in Proposition \ref{P1}, $\Lambda^{\prime}(0, 0, 0)$ is surjective. We can therefore apply Lyusternik-Graves' Theorem \ref{Lyusternik} to this data, in particular, for any initial data $(z_{0},k_{0}) \in H^{1}_{0}(-1,1)\times\mathbb{R}$ with sufficiently small norm
\[
\|(z_{0},k_{0})\|_{H^{1}_{0}(-1,1)\times\mathbb{R}} < \varepsilon,
\]
for some sufficiently small $\varepsilon>0$, there exists $(z,k,v)\in\mathds{X}$ such that
\[
\Lambda(z,k,v) = (0,0,z_{0},k_{0}).
\] 
Consequently, $(z,k)$ is the solution of \eqref{s6} associated to the control $v$, and 
the exponential growth of the weight $\rho_{0}$ as $t\rightarrow T^{-}$ and $\rho_{0}(z,k)\in L^{2}(\tilde{Q})\times L^{2}(0,T)$, ensures that
\begin{eqnarray*}
	z(\cdot,T)=0\quad\mbox{and}\quad k(T)=0.
\end{eqnarray*}
\textcolor{red}{
Moreover, we have
\begin{align*}
	\|(z,k,v)\|_{\mathds{X}}\le C (\|(z_0, k_0)\|_{H^1(-1,1)\times \mathbb{R}},
	\end{align*}
for some constant $C>0$ independent of $\varepsilon$. }
This leads to the proof of Theorem \ref{result} and, as explained in the introduction and Remark \ref{reg of contr}, we immediately obtain the main Theorem \ref{main result}.

%%%%%%%%%%%%%%%%%%%%%%%%%%%%%%%%%%%%%%%%%%%%%%%
\section{Numerical Modeling with PINNs
	%Numerics
} \label{Sec 4}
%%%%%%%%%%%%%%%%%%%%%%%%%%%%%%%%%%%%%%%%%%%%%%%
In recent years, Physics-Informed-Neural-Networks (PINNs) have attracted considerable attention for solving mathematical models based on nonlinear PDEs; see, for example, \cite{PENG2022106067, ZHANG2024108229} and the references therein. This section addresses the numerical resolution of problem \eqref{s1}, and its linearized counterpart, by using Physics-Informed-Neural-Networks (PINNs). 

%%%%%%%%%%%%%%%%%%%%%%%%%%%%%%%%%%%%%%%%%%%%%%%
\subsection{Description of the numerical algorithm}
%%%%%%%%%%%%%%%%%%%%%%%%%%%%%%%%%%%%%%%%%%%%%%% 

Fostered by the seminal paper \cite{rai19},  PINNs algorithm was adapted to deal with controllability problems involving partial differential equations in \cite{per23}.  The underlying idea is rather simple: it amounts to use neural networks as surrogates of the state and control variables, and then minimizing the residual of the ODE-PDE system as well as its boundary, initial and terminal conditions. Next, we provide more specific details. 

For the sake of clarity in the exposition, let us consider the linearized system \eqref{s5}, with $\overline{\psi}=0, \overline{h}=1$ (which corresponds to Experiment 1 below), the nonlinear case being completely analogous. More precisely, we consider the system  

\begin{equation*}
	\left\{
	\begin{aligned}
		&z_{t}-z_{xx} -z= \mathds{1}_{\omega}v & & \text {in}\; \Omega_T := \Omega \times (0, T), \\
		&z=0  & & \text {on}\; \Gamma_T :=\partial\Omega\times (0,T), \\
		&k^{\prime}(t)+2\mu z_{x}(1,t)=0 & & \text {in}\;(0,T),  \\
		& k(0)=k_{0},\\
		& z(\cdot,0)=z_{0} & &\text {in}\; \Omega,
	\end{aligned}
	\right.
	%\label{model_linear}
\end{equation*}
where $\Omega\subset \mathbb{R}^{d}$ ($d=1$) is the spatial domain ($\Omega = (-1, 1)$, in Experiment 1 below).

The PINNs algorithm is structured as follows:

\textbf{Step 1: Neural networks. } Among different possibilities, we consider deep feedforward neural networks (also known in the literature as Multilayer Perceptrons (MLPs))  with $d+1$ input channels ($\boldsymbol{x} =(x,t)\in \mathbb{R}^{d+1}$ for the surrogate of $z(x, t)$ and $v(x, t)$, and just $ t$ for the input variable of $k(t)$,  and an scalar output (say $\hat{z}$ to fix ideas). More specifically, $\hat{z}\left(x,t;\theta\right)$  is constructed as
\begin{equation*}
	%\label{NN}
	\left\{\begin{array}{ll}
		\text{input layer:} & \mathcal{N}^0(x)=\bs x=(x,t)\in \mathbb{R}^{d+1}, \\
		\text{hidden layers:} & \mathcal{N}^{\ell}(\bs{x})=\sigma\left(\bs{W}^{\ell}\mathcal{N}^{\ell -1}(\bs{x})+\bs{b}^{\ell}\right)\in\mathbb{R}^{N_{\ell}}, \quad \ell = 1, \cdots, \mathfrak{L}-1, \\
		\text{output layer:} & \hat{z}\left(\bs{x};\bs{\theta}\right) = \mathcal{N}^\mathfrak{L}(\bs{x})=\bs{W}^{\mathfrak{L}}\mathcal{N}^{\mathfrak{\mathfrak{L}} -1}(\bs{x})+\bs{b}^{\mathfrak{L}}\in\mathbb{R},
	\end{array}	
	\right.
\end{equation*}
where 
\begin{itemize}
	\item $\mathcal{N}^{\ell}:\mathbb{R}^{d_{in}}\rightarrow\mathbb{R}^{d_{out}}$ is the $\ell$ layer with $N_{\ell}$ neurons, 
	
	\item  $\bs{W}^{\ell}\in\mathbb{R}^{N_{\ell}\times N_{\ell -1}}$ and
	$\bs{b}^{\ell}\in \mathbb{R}^{N_{\ell}}$ are, respectively, the weights and biases so that $\bs{\theta}=\left\{ \bs{W}^{\ell}, \bs{b}^{\ell}\right\}_{1\leq\ell\leq \mathfrak{L}}$ are the parameters of the neural network, and  
	\item  $\sigma$ is an activation function, which acts component-wise. Throughout this paper, we consider smooth activation functions such as hyperbolic tangent $\sigma (s)=\tanh (s)$, with $s\in\mathbb{R}$. 	
\end{itemize}

The reason we use this neural network architecture is that MLPs have been proven to be universal approximators for smooth functions \cite{pinkus1999approximation}. 

\textbf{Step 2: Training dataset. } A dataset $\mathcal{T}$ of scattered points is selected in the interior domains $\mathcal{T}_{\text{int}, \Omega_T}\subset \Omega_T$ and $\mathcal{T}_{\text{int}, T}\subset (0, T)$, and on the boundaries $\mathcal{T}_{\Gamma}\subset \Gamma_T$, $\mathcal{T}_{t=0} \subset \Omega\times \left\{ 0\right\}$, $\mathcal{T}_{t=T}\subset \Omega\times \left\{ T\right\}$. Thus, $\mathcal{T}=\mathcal{T}_{\text{int}, \Omega_T}\cup\mathcal{T}_{\text{int}, T}\cup\mathcal{T}_{\Gamma}\cup\mathcal{T}_{t=0}\cup\mathcal{T}_{t=T}$. The number of selected points in $\mathcal{T}_{\text{int}, \Omega_T}$ and  $\mathcal{T}_{\text{int}, T}$ are, respectively, denoted by $N_{int, \Omega_T}$ and $N_{int, T}$. Analogously, $N_{b}$ is the number of points on the boundary $\Gamma$, and  $N_{0}$ and $N_{T}$ stand for the number of points in $\mathcal{T}_{t=0}$ and $\mathcal{T}_{t=T}$, respectively. The total number of collocation nodes is denoted by $N$ and we write $\mathcal{T}_N$ instead of $\mathcal{T}$ to indicate more dearly the number of points $N$ used hereafter. 

\textbf{Step 3: Loss function. }  A weighted summation of the $L^2$ norm of residuals for the involved differential equations, boundary, initial and final conditions is considered as the loss function to be minimized during the training process. It is composed of the following eight terms: given a neural networks approximations $\hat{z}$, $\hat{k}$ and $\hat{v}$ of $z$, $k$ and $v$, respectively, we define:

\[
%\scriptsize
%\small
\begin{array}{lll}
	\mathcal{L}_{\text{int}, PDE}\left( \bs{\theta};\mathcal{T}_{\text{int, } \Omega_T}\right) & = \displaystyle\sum_{j=1}^{N_{\text{int, } \Omega_T}} w_{j,\text{int, } \Omega_T}\vert \hat{z}_{t}(\bs{x}_j)-\hat{z}_{xx}(\bs{x}_j) - \hat{z}(\bs{x}_j)- \mathds{1}_{\omega}\hat{v}(\bs{x}_j) \vert^2 , & \bs{x}_j\in \mathcal{T}_{\text{int, } \Omega_T}, \\
	\\
	\mathcal{L}_{\text{int }, ODE}\left( \bs{\theta};\mathcal{T}_{\text{int, } T}\right) & = \displaystyle\sum_{j=1}^{N_{\text{int, } T}} w_{j,\text{int, } T}\vert \hat{k}^{\prime}(\bs{x}_j)+2\mu \hat{z}_{x}(1, \bs{\tilde{x}}_j) \vert^2 , & (\bs{x}_j, \bs{\tilde{x}}_j) \in \mathcal{T}_{\text{int, } T}\times \mathcal{T}_{\text{int, } \Omega_T}, \\
	\\
	\mathcal{L}_{\Gamma}\left( \bs{\theta};\mathcal{T}_{\Gamma}\right) & = \displaystyle\sum_{j=1}^{N_{b}} w_{j, b}\vert \hat{z}(\bs{x}_j;\bs{\theta}) \vert^2, & \bs{x}_j\in \mathcal{T}_{\Gamma}, \\
	\\
	\mathcal{L}_{t=0}^z\left( \bs{\theta};\mathcal{T}_{t=0}\right) & = \displaystyle\sum_{j=1}^{N_{0}} w_{j, 0}\vert \hat{z}(\bs{x}_j;\bs{\theta}) - z_0(\bs{x}_j) \vert^2, &  \bs{x}_j\in \mathcal{T}_{t=0}, \\
	\\
	\mathcal{L}_{t=0}^k\left( \bs{\theta}\right) & = \vert \hat{k}(0; \bs{\theta}) - k_0\vert, \\
	\\	
	\mathcal{L}_{t=T}^z\left( \bs{\theta};\mathcal{T}_{t=T}\right) & = \displaystyle\sum_{j=1}^{N_{T}} w_{j, T}\vert \hat{z}(\bs{x}_j;\bs{\theta})  \vert^2, &  \bs{x}_j\in \mathcal{T}_{t=T},  \\
	\\
	\mathcal{L}_{t=T}^{k}\left( \bs{\theta}\right) & = \vert \hat{k}(T;\bs{\theta})  \vert, &  
	\\
	\\
	\mathcal{L}_{\hat{z}\geq 0}(\bs{\theta})	& = \displaystyle\sum_{j=1}^{N_{\text{int, } \Omega_T}} w_{j,\text{int, } \Omega_T} \frac{1}{\beta}\log\left(1 +  e^{-\beta\hat{z}(\bs{x}_j)}\right)\vert^2, & \bs{x}_j\in \mathcal{T}_{\text{int, } \Omega_T},
\end{array}
\]	
where the $ w_{s}$ are the weights corresponding to selected quadrature rules.

\begin{figure}[h!]
	\centering
	\boxed{
		\begin{tikzpicture}[
			every node/.style={draw, thick, minimum width=1.2cm, minimum height=1cm, align=center, rounded corners},
			node distance=1.6cm and 2.6cm, 
			col1/.style={fill=blue!20},
			col2/.style={fill=green!20},
			col3/.style={fill=orange!20},
			col4/.style={fill=red!20},
			col5/.style={fill=purple!20},
			arrow/.style={thick, ->, >=stealth}
			]
			
			% Colonnes avec des coordonnées spécifiques
			\node[col1] (C11) at (-1, 0.9) {\(\hat{z}(x,t)\)};
			\node[col1] (C12) at (-1, -0.4) {\(\hat{k}(t)\)};
			\node[col1] (C13) at (-1, -1.6) {\(\hat{v}(x,t)\)};
			
			\node[col2] (C21) at (4, 2) {\(\hat{z}_{t}-\hat{z}_{xx} -\hat{z}-\mathds{1}_{\omega}\hat{v} =0 \)\\ 
				\(\hat{k}^{\prime}(t)+2\mu \hat{z}_{x}(1,t)=0 \)};
			\node[col2] (C22) at (4, 0.5) {\(\hat{z}(-1,t)=\hat{z}(1,t)=0  \)};
			\node[col2] (C23) at (4, -1) {\(\hat{z}(x,0)=\hat{z}_{0}(x)\)\\ \(\hat{k}(0)=\hat{k}_0\)};
			\node[col2] (C24) at (4, -2.5) {\(\hat{z}(x,T)=\hat{z}_{T}(x)\)\\ \(\hat{k}(T)=\hat{k}_T\)};
			
			\node[col3] (C31) at (8.5, -0.2) {\(\displaystyle\min_{\theta}\mathcal{L}(\hat{z},\hat{k},\hat{v};\theta)\)};
			
			\node[col4] (C41) at (11.2, -0.2) {\(\theta^{\ast}\)};
			
			\node[col5] (C51) at (13.5, 1.3) {\(\hat{z}(x,t; \theta^{\ast})\)};
			\node[col5] (C52) at (13.5, -0.2) {\(\hat{k}(t, \theta^{\ast})\)};
			\node[col5] (C53) at (13.5, -1.7) {\(\hat{v}(x,t; \theta^{\ast})\)};
			
			%			\node[col6] (C61) at (15, -1.57) {C6.1};
			
			% Connexions entre les n?uds
			\draw[arrow, draw=gray, thin] (C11) -- (C21);
			\draw[arrow, draw=gray, thin] (C11) -- (C22);
			\draw[arrow, draw=gray, thin] (C11) -- (C23);
			\draw[arrow, draw=gray, thin] (C11) to[out=-40, in=180] (C24);
			\draw[arrow, draw=gray, thin] (C12) -- (C21);
			%			\draw[arrow, draw=gray, thin] (C12) -- (C22);
			\draw[arrow, draw=gray, thin] (C12) -- (C23);
			\draw[arrow, draw=gray, thin] (C12) -- (C24);
			
			\draw[arrow, draw=gray, thin] (C13) to[out=40, in=180] (C21);
			
			\draw[arrow, draw=gray, thin] (C21) -- (C31);
			\draw[arrow, draw=gray, thin] (C22) -- (C31);
			\draw[arrow, draw=gray, thin] (C23) -- (C31);
			\draw[arrow, draw=gray, thin] (C24) -- (C31);
			
			\draw[arrow, draw=gray, thin] (C31) -- (C41);
			
			\draw[arrow, draw=gray, thin] (C41) -- (C51);
			\draw[arrow, draw=gray, thin] (C41) -- (C52);
			\draw[arrow, draw=gray, thin] (C41) -- (C53);

	\end{tikzpicture} }
	\caption{PINN-based algorithm for approximating the exact states and control. The neural networks $\hat{z}\left(x,t,\bs{\theta}\right)$, $\hat{k}\left(t,\bs{\theta}\right)$ and $\hat{v}\left(x,t,\bs{\theta}\right)$ are required to satisfy, in the least squares sense, the ODE-PDE system, boundary condition, initial conditions, and exact controllability conditions. Then,  the residual on training points $\mathcal{L}\left(\bs{\theta};\mathcal{T}\right)$ is minimized to get the optimal set of parameters $\bs{\theta}^\ast$ of the neural networks. This leads to the PINN exact states $\hat{z}\left(x,t;\bs{\theta}^\ast\right)$, $\hat{k}\left(t;\bs{\theta}^\ast\right)$ as well as the PINN control $\hat{v}\left(x,t;\bs{\theta}^\ast\right)$.}
	\label{fig:aligned_diagram}
\end{figure}
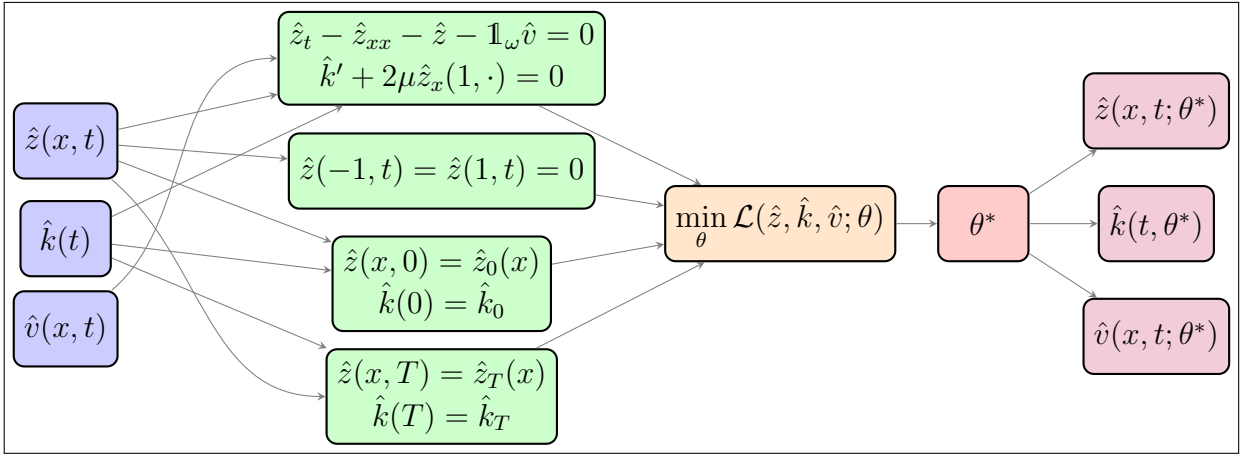

The loss function used for training is 
\begin{equation}
	\label{loss_total}
	\begin{array}{ll}
		\mathcal{L}\left( \bs{\theta};\mathcal{T}\right) & = 
		w_1\mathcal{L}_{\text{int}, PDE}\left( \bs{\theta};\mathcal{T}_{\text{int, } \Omega_T}\right) + w_2\mathcal{L}_{\text{int }, ODE}\left( \bs{\theta};\mathcal{T}_{\text{int, } T}\right)\\
		& + w_3\mathcal{L}_{\Gamma}\left( \bs{\theta};\mathcal{T}_{\Gamma}\right) \\
		& + w_4\mathcal{L}_{t=0}^z\left( \bs{\theta};\mathcal{T}_{t=0}\right)
		+ w_5\mathcal{L}_{t=0}^k\left( \bs{\theta}\right) \\
		& + w_6\mathcal{L}_{t=T}^z\left( \bs{\theta};\mathcal{T}_{t=T}\right) 
		+ w_7\mathcal{L}_{t=T}^{k}\left( \bs{\theta}\right)\\
		& +	w_8\mathcal{L}_{\hat{z}\geq 0}(\bs{\theta}),
	\end{array} 
\end{equation}
where, as above, the $w_s$ are suitable positive weights. See Figure \ref{fig:aligned_diagram} for a schematic illustration.

Crucially, the definitions of $\mathcal{L}_{t=0}^k\left( \bs{\theta}\right)$, $\mathcal{L}_{t=T}^{k}\left( \bs{\theta}\right)$, and $\mathcal{L}_{\hat{z}\geq 0}(\bs{\theta})$ are rooted by the fact that we are interested in having good approximations of the original variables $L(t) = k(t)^{1/2}$, and $\varrho(x, t)$. In particular, to enforce non-negativity of the density $\varrho(x, t)$, the log-sum-exp approximation of the max function (i.e, $\max (x_1, x_2) \approx \frac{1}{\beta} \log\left(e^{\beta x_1} + e^{\beta x_2}\right)$, with $\beta >>1$) has been considered.

\textbf{Step 4: Training process. } The final step in the PINN algorithm amounts to minimizing \eqref{loss_total}, i.e.,
\begin{equation*}
	%\label{train_Fisher}
	\bs{\theta}^\ast = \arg\min_{\bs{\theta}}	\mathcal{L}\left( \bs{\theta};\mathcal{T}\right).
\end{equation*}

%%%%%%%%%%%%%%%%%%%%%%%%%%%%%%%%%%%%%%%%%%%%%%%
\subsection{Numerical experiments}
%%%%%%%%%%%%%%%%%%%%%%%%%%%%%%%%%%%%%%%%%%%%%%%

This section provides numerical simulation results for both the nonlinear Fisher-Stephan problem \eqref{s0} and its linearized approximation (system \eqref{s5}) . 

In all experiments that follow, a multi-layer neural network, as described in the preceding section, with the $\tanh$ as activation function, is used to approximate $z$, $v$, and $k$. In the three cases, the neural network is composed of two hidden layers with fifty neurons in each layer.  
Uniformly distributed random quadrature nodes are employed for training the neural network. For simplicity, a simple Monte Carlo integration has been implemented. The training process, i.e. minimization of $\mathcal{L}\left(\bs{\theta};\mathcal{T}_N\right)$, is carried out with the  ADAM optimizer \cite{kingma2014adam} with learning rate $10^{-3}$. The required gradients are computed by using Automatic Differentiation. By default, a Glorot uniform distribution \cite{glorot2010understanding} is used as initial weights. The biases are  initialized  with zero value.
%\par
%In the following experiments, we aim to illustrate the main theorem \ref{main result} using the trajectories:
%\begin{eqnarray*}
%	\overline{\varrho}=0\quad\mbox{and}\quad \overline{L}=1.
%\end{eqnarray*}
%We aim to find a control $u$

%%%%%%%%%%%%%%%%%%%%%%%%%%%%%%%%%%%%%%%%%%%%%%%
\subsubsection{Experiment 1: Linearized system \eqref{s5}}
%%%%%%%%%%%%%%%%%%%%%%%%%%%%%%%%%%%%%%%%%%%%%%%
This section shows numerical simulation results for Theorem \ref{main result} using linearized system \eqref{s5} and for the following data: $\Omega=(-1, 1)$, $\overline{\varrho}=0, \overline{L}=1$, ( then $\overline{\psi}=0, \overline{h}=1$, $m=0, n=-1$ ), $\mu = 0.5$, $T=1$ and $f_0 = g_0 = 0$. The initial conditions are  $k_0=-0.5$ (so that $L_0 = \sqrt{0.5})$, $z_0=\left\{ \begin{array}{ll}
	\sin (\pi x), & 0\leq x\leq 1, \\
	0,	& -1\leq x < 0,
\end{array}\right. $
and the control region is $\omega = (-1 +\varepsilon, -\varepsilon )$, with $\varepsilon=0.1$.

Table \ref{table:exp1a_regu} collects results for the evolution, with respect to the number of training points, of all the losses in \eqref{loss_total} as well as generalization error (difference between training and test errors). $20000$ epochs have been performed in all cases. Since our domain is rectangular, training points are the result of the Cartesian product of points in the spatial and time directions. The same number of random points has been taken in both directions. It is observed that beyond a certain number of training points (on the order of $10^6$), further increasing the number of points does not significantly enhance the numerical results.

\begin{table*}[h!]
	\begin{center}
		\caption{ Experiment 1: Training errors for the different contributions in \eqref{loss_total}. Last column corresponds to generalization error ($\mathcal{E}_{\text{gener}}$).    }
		\label{table:exp1a_regu}
		\tiny
		\begin{tabular}{lccccccccc}
			\hline\noalign{\smallskip}
			&  $\mathcal{L}_{\text{int}, PDE}$ &   $ \mathcal{L}_{\text{int }, ODE}$ & $\mathcal{L}_{\Gamma}$  & $\mathcal{L}_{t=0}^z$ & $ \mathcal{L}_{t=0}^k$ & $\mathcal{L}_{t=T}^z$ & $\mathcal{L}_{t=T}^{k}$ & $\mathcal{L}_{\hat{z}\geq 0}$ & $\mathcal{E}_{\text{gener}}$ \\
			\\
			$ N = 10^6$ & $9.0\times 10^{-4}$ & $1.0\times 10^{-4}$  & $4.0\times 10^{-4}$ & $3.0\times 10^{-4}$ & $1.5\times 10^{-3}$ & $4.0\times 10^{-4}$ & $3.3 \times 10^{-3}$ & $4.0 \times 10^{-2}$ & $1.5\times 10^{-2}$\\ 
			$ N = 10^8 $ &  $6.0\times 10^{-4}$ & $1.0\times 10^{-5}$ & $1.0\times 10^{-5}$ & $2.0\times 10^{-4}$ & $1.6\times 10^{-3}$ & $1.0\times 10^{-5}$ & $3.9\times 10^{-3}$ & $4.5\times 10^{-2}$& $3.4\times 10^{-3}$ \\ 
			\noalign{\smallskip}\hline
		\end{tabular}
	\end{center}
\end{table*}

\begin{figure}[H]
	\centering
	\begin{tabular}{cc}
		\hspace{-1cm}\includegraphics[height=70mm]{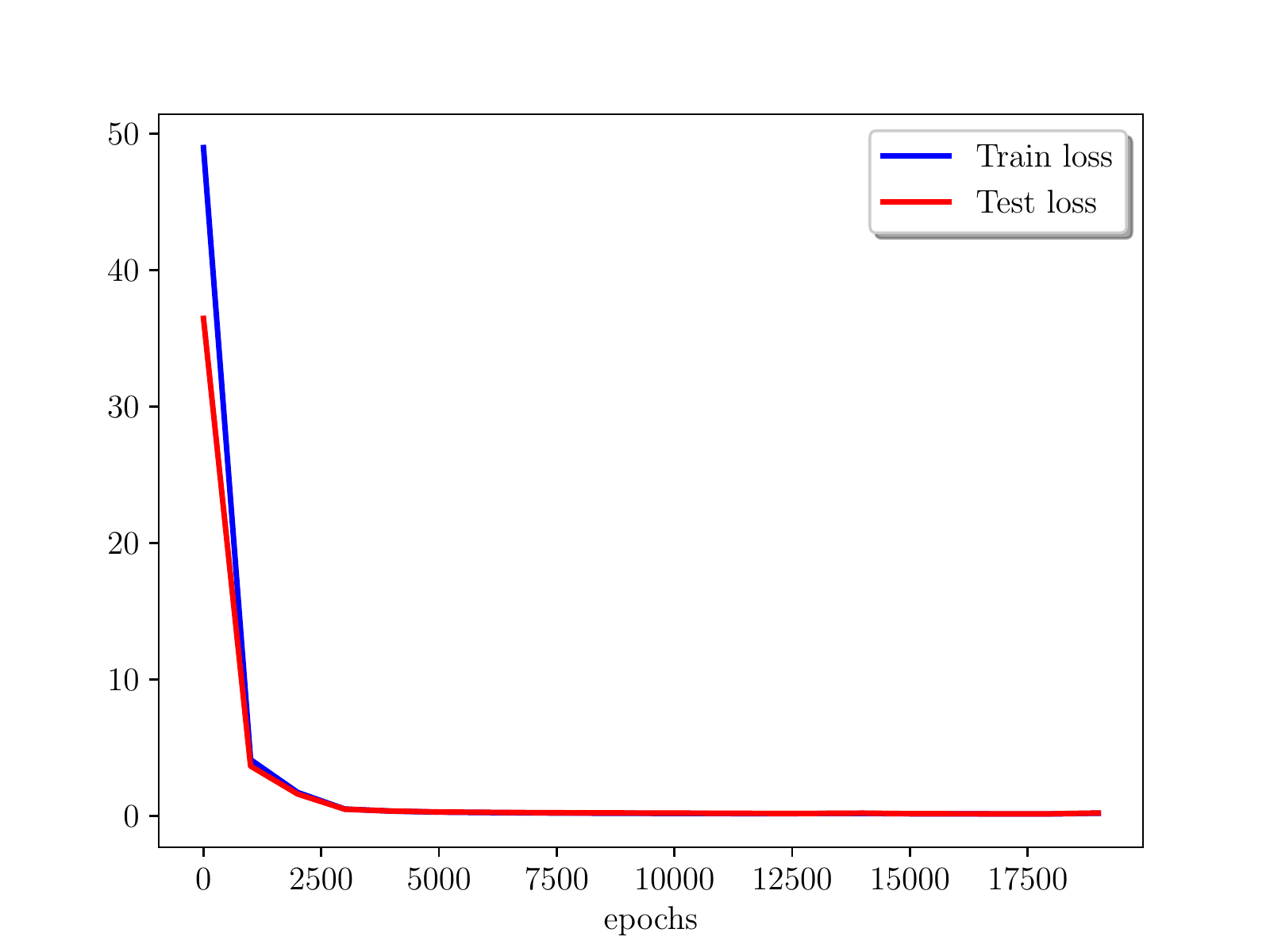}  & \hspace{-1cm}\includegraphics[height=70mm]{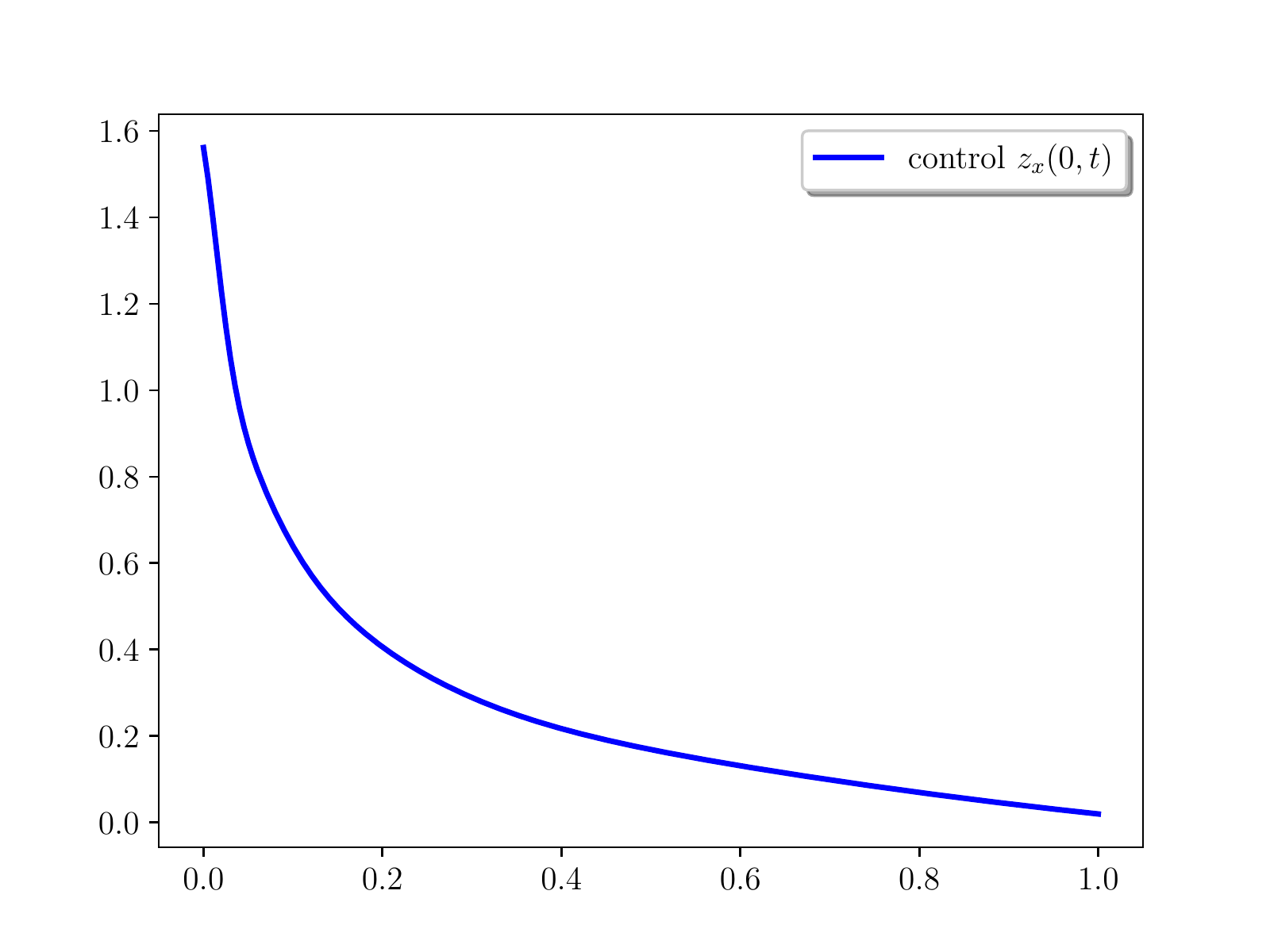}
	\end{tabular}	
	\caption{Experiment 1. Convergence history \textbf{(left)}, and control \textbf{(right)}.  }
	\label{exp1:history_control}
\end{figure}

Figure \ref{exp1:history_control} (left) depicts the evolution of train and test losses. It is observed that neither under-fitting nor over-fitting take place. The control variable, which is computed from the PINN algo desribed in the preceding section, is shown in \ref{exp1:history_control} (right). 

%\begin{figure}[H]
%	\centering
%		\includegraphics[height=70mm]{linear_converge_history.png}  
%	\caption{Experiment 1. $T=1$. Convergence history.}
%	\label{exp1:history}
%\end{figure}
%Notice that system \eqref{s5} is decoupled, meaning that we can solve first for $z$ and then for $k$. Consequently, it is observed that the influence of the value of $n$ is negligible. In both cases, the system is dissipative w.r.t. $z$.   

\begin{figure}[H]
	\centering
	\begin{tabular}{cc}
		\hspace{-1cm}\includegraphics[height=70mm]{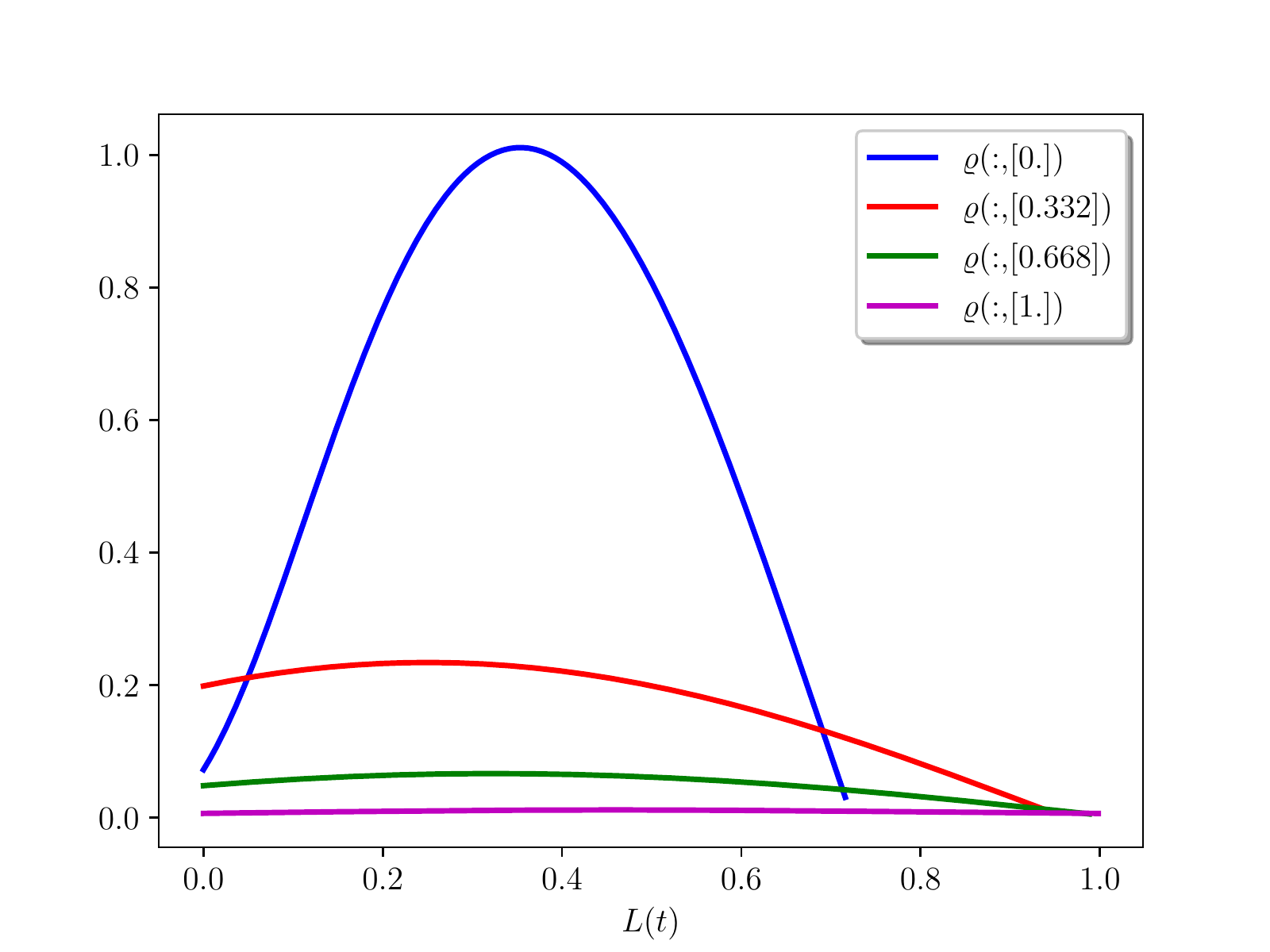}  & \hspace{-1cm}\includegraphics[height=70mm]{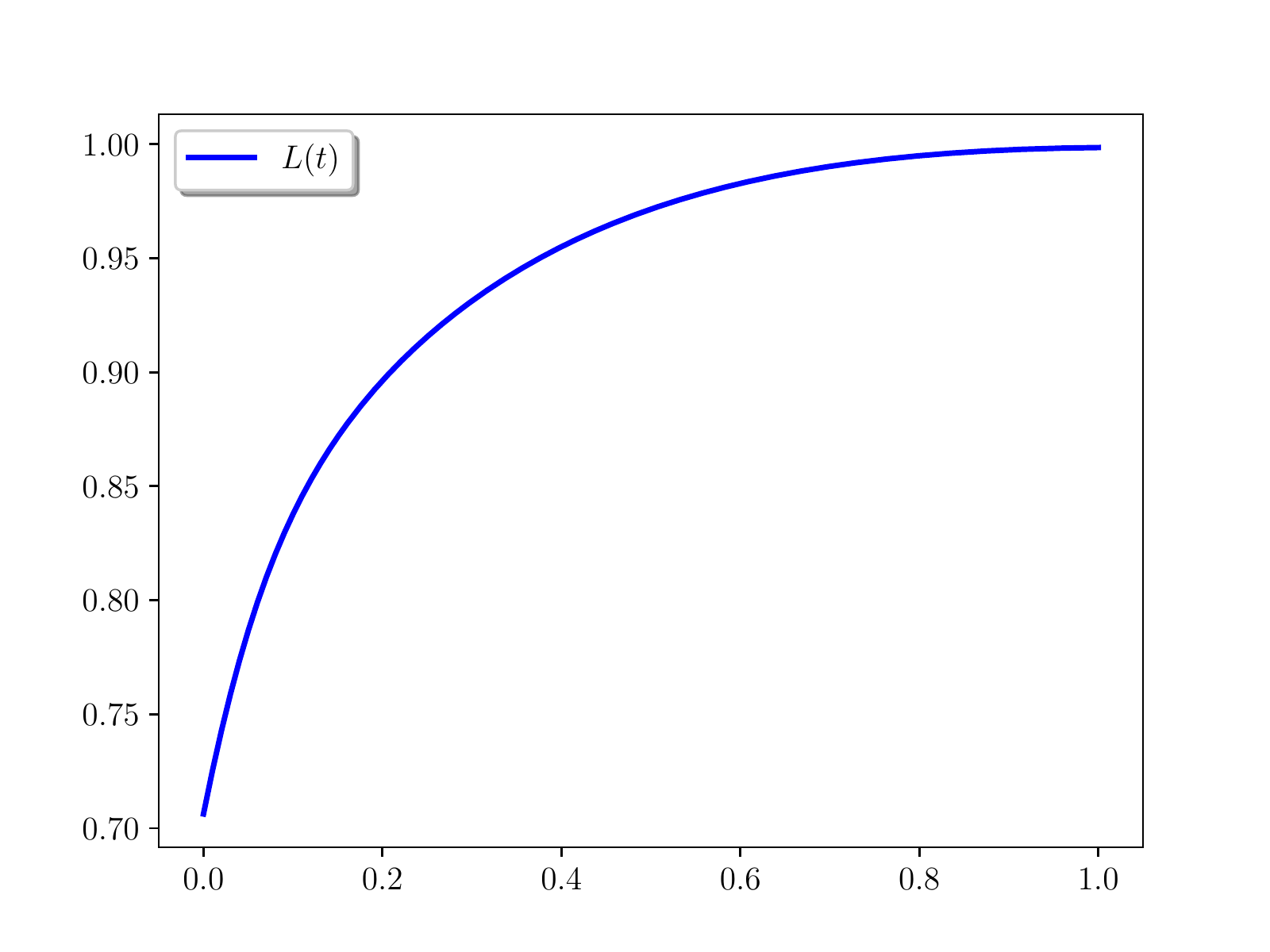} 
	\end{tabular}	
	\caption{Experiment 1.  \textbf{(Left)} Plot of the densities associated with the linearized system at different times, and  state variable $L(t)$ \textbf{(right)}.  }
	\label{exp1:linear_states}
\end{figure}

Finally, Figure \ref{exp1:linear_states} (left) displays results for the controlled density $\varrho(xL(t), t)$ at several discrete times (left). Importantly, all the densities remain non-negative as a consequence of the penalization term $\mathcal{L}_{\hat{z}\geq 0}$ introduced in the loss function \eqref{loss_total}. The  state variable $L(t)$ is plotted in Figure \ref{exp1:linear_states} (right). In both cases, it is observed that all the considered initial and controllability conditions  are satisfied with accuracy.

%%%%%%%%%%%%%%%%%%%%%%%%%%%%%%%%%%%%%%%%%%%%%%%%%%%%%%%%%%%%%%%%
\subsubsection{Experiment 2: Nonlinear system \eqref{s1}}
%%%%%%%%%%%%%%%%%%%%%%%%%%%%%%%%%%%%%%%%%%%%%%%%%%%%%%%%%%%%%%%

This section reviews on numerical simulation results for system \eqref{s1}, which are obtained from the PINN approximation of \eqref{s2}. For comparison with its associated linearized system, we take the same initial and final conditions, namely, $h_0=0.5$, $\mu = 0.5$, $T=1$, and $\psi_0 (x)=
\sin (\pi x)$.

\begin{figure}[H]
	\centering
	\begin{tabular}{cc}
		\hspace{-1cm}\includegraphics[height=70mm]{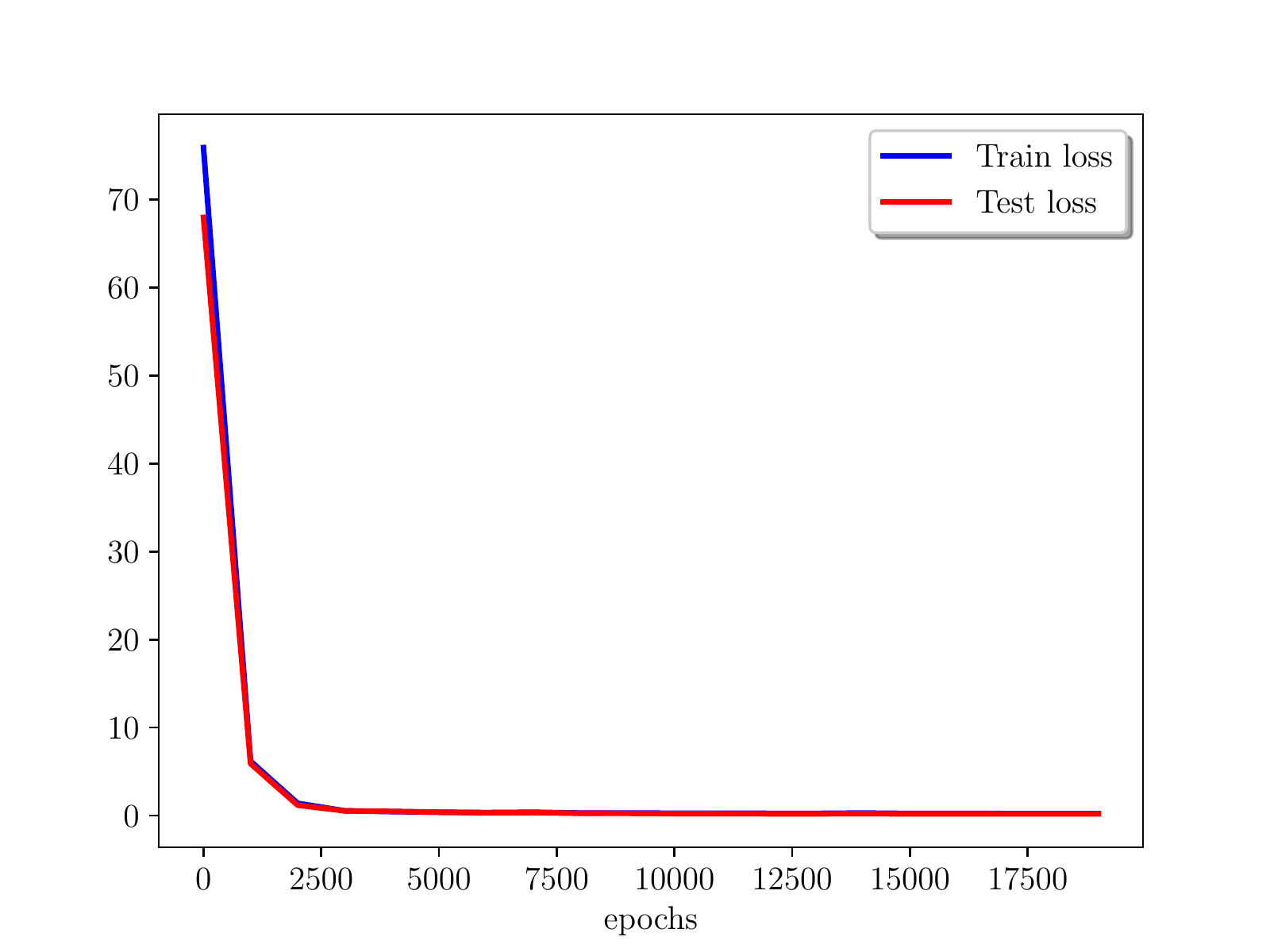}  & \hspace{-1cm}\includegraphics[height=70mm]{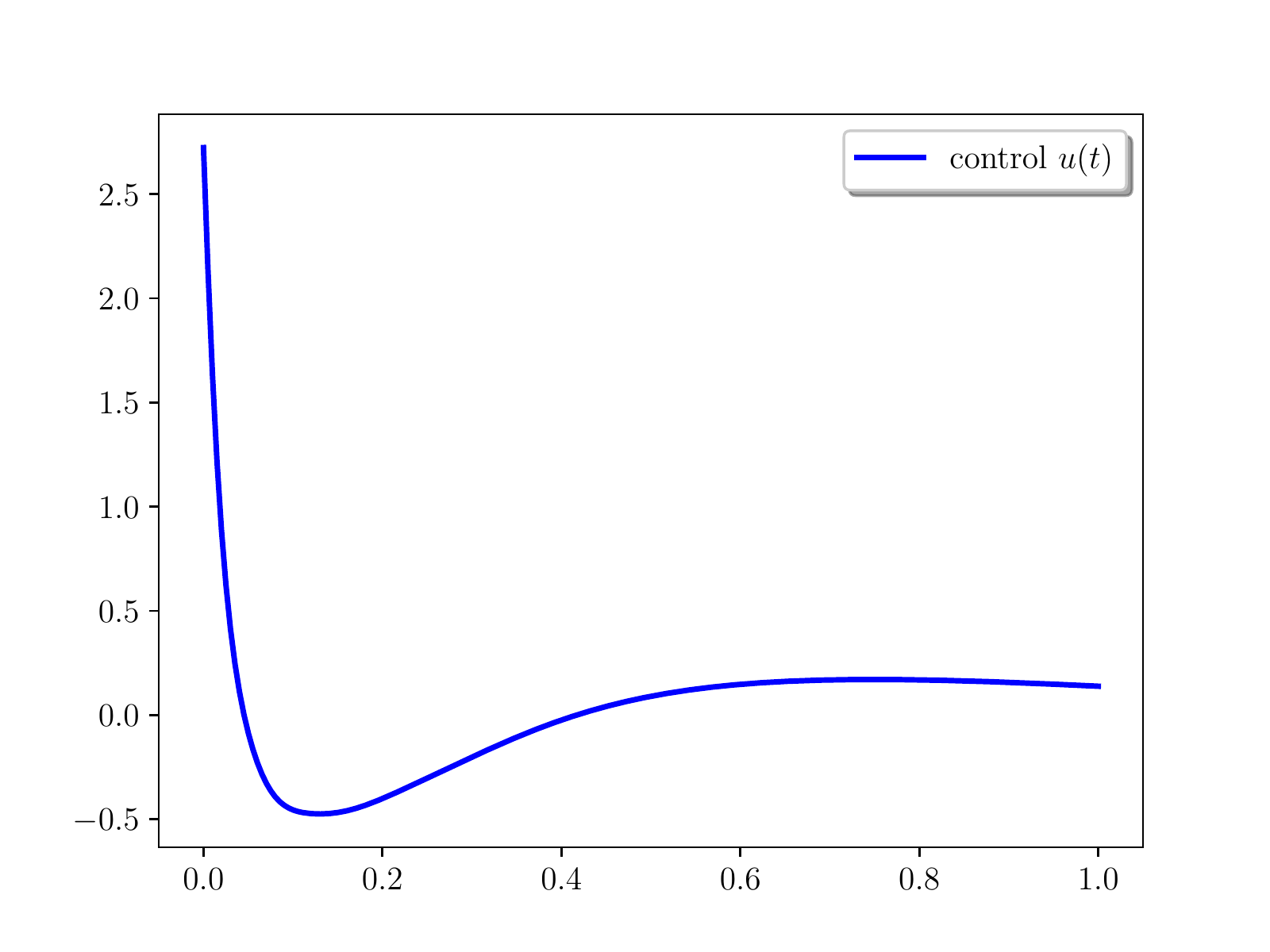}
	\end{tabular}	
	\caption{Experiment 2. $T=1$. \textbf{(Left)} Convergence history. \textbf{(Right)} Control $u(t)$.  }
	\label{exp2:history_and_control}
\end{figure}

As in the linearized case, the simulation results that follow have been obtained with $10^6$ random training points in both the spatial and time dimensions. As for test, $10^2$ random points are used. Figure \ref{exp2:history_and_control} (left) shows the convergence history of the PINNs algorithm. Table \ref{table:exp_2} collects numerical results for the different contributions of the loss function after convergence of the algorithm for $10^6$ and $10^8$. 

\begin{table*}[h!]
	\begin{center}
		\caption{ Experiment 2: Training errors for the different contributions in the loss function. Last column corresponds to generalization error ($\mathcal{E}_{\text{gener}}$).    }
		\label{table:exp_2}
		\tiny
		\begin{tabular}{lccccccccc}
			\hline\noalign{\smallskip}
			&  $\mathcal{L}_{\text{int}, PDE}$ &   $ \mathcal{L}_{\text{int }, ODE}$ & $\mathcal{L}_{\Gamma}$  & $\mathcal{L}_{t=0}^{\psi}$ & $ \mathcal{L}_{t=0}^h$ & $\mathcal{L}_{t=T}^{\psi}$ & $\mathcal{L}_{t=T}^{h}$ & $\mathcal{L}_{\hat{z}\geq 0}$ & $\mathcal{E}_{\text{gener}}$ \\
			\\
			$ N = 10^6$ & $6.6\times 10^{-3}$ & $4.0\times 10^{-4}$  & $3.0\times 10^{-3}$ & $1.0\times 10^{-4}$ & $1.4\times 10^{-3}$ & $1.0\times 10^{-4}$ & $3.3 \times 10^{-3}$ & $2.7 \times 10^{-2}$ & $4.2\times 10^{-2}$\\ 
			$ N = 10^8 $ &  $7.4\times 10^{-3}$ & $4.0\times 10^{-4}$ & $3.0\times 10^{-4}$ & $1.0\times 10^{-4}$ & $1.1\times 10^{-3}$ & $1.0\times 10^{-4}$ & $1.8\times 10^{-3}$ & $2.6\times 10^{-2}$& $7.5\times 10^{-3}$ \\ 
			\noalign{\smallskip}\hline
		\end{tabular}
	\end{center}
\end{table*}

Figure \ref{exp2:history_and_control} (right) depicts the computed control variable.

%\begin{figure}[H]
%	\centering
%	\includegraphics[height=70mm]{nonlinear_convergence_history.pdf}  
%	\caption{Experiment 2. $T=1$. Convergence history.}
%	\label{exp2:history}
%\end{figure}

%Similarly to the linearized case, dissipation w.r.t. $\psi$ is also observed. However, without control, $h$ increases, as shown in Figure \ref{exp2:state_L} (left).
\begin{figure}[H]
	\centering
	\begin{tabular}{cc}
		\hspace{-1cm}\includegraphics[height=70mm]{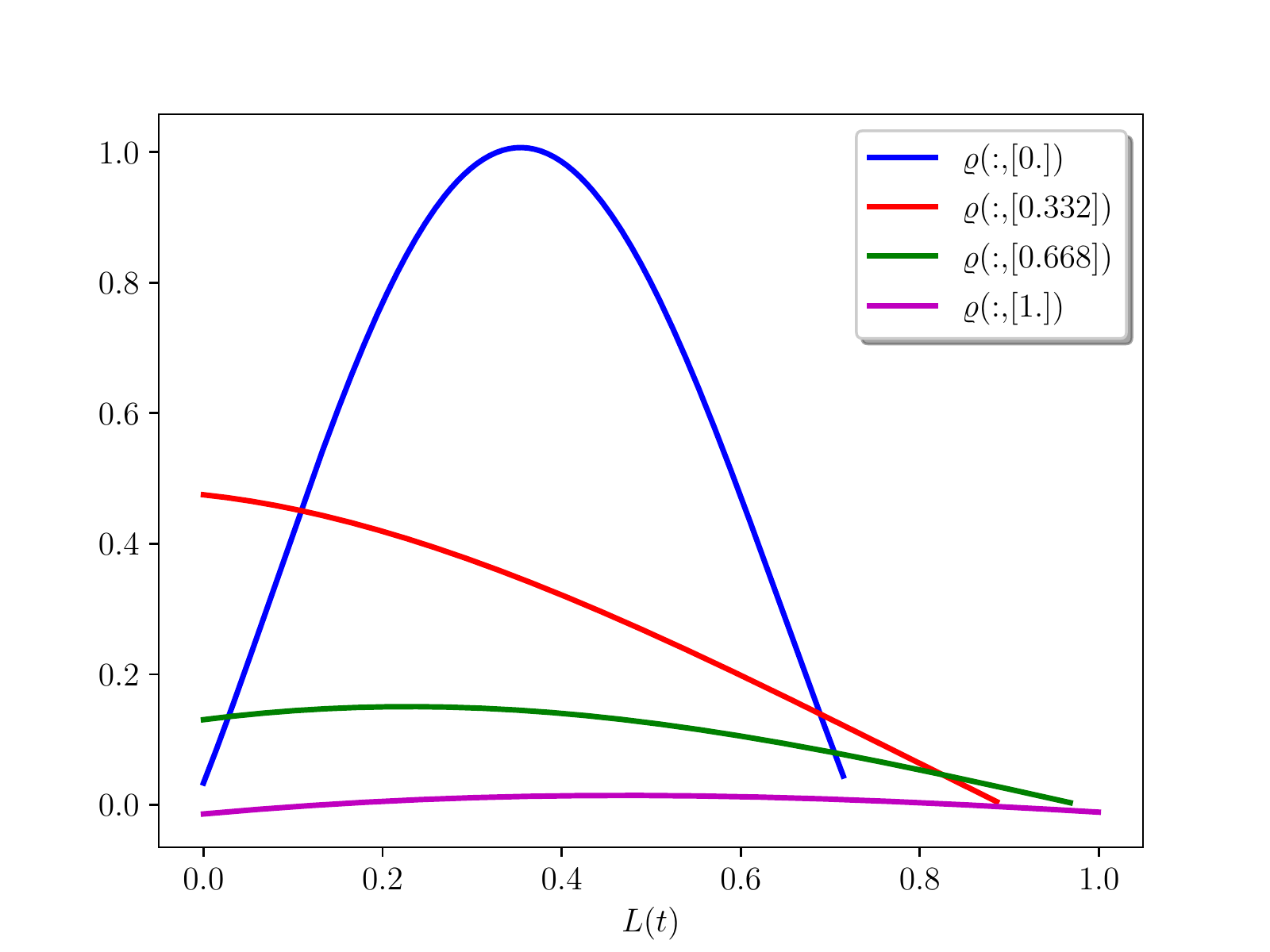}  & \hspace{-1cm}\includegraphics[height=70mm]{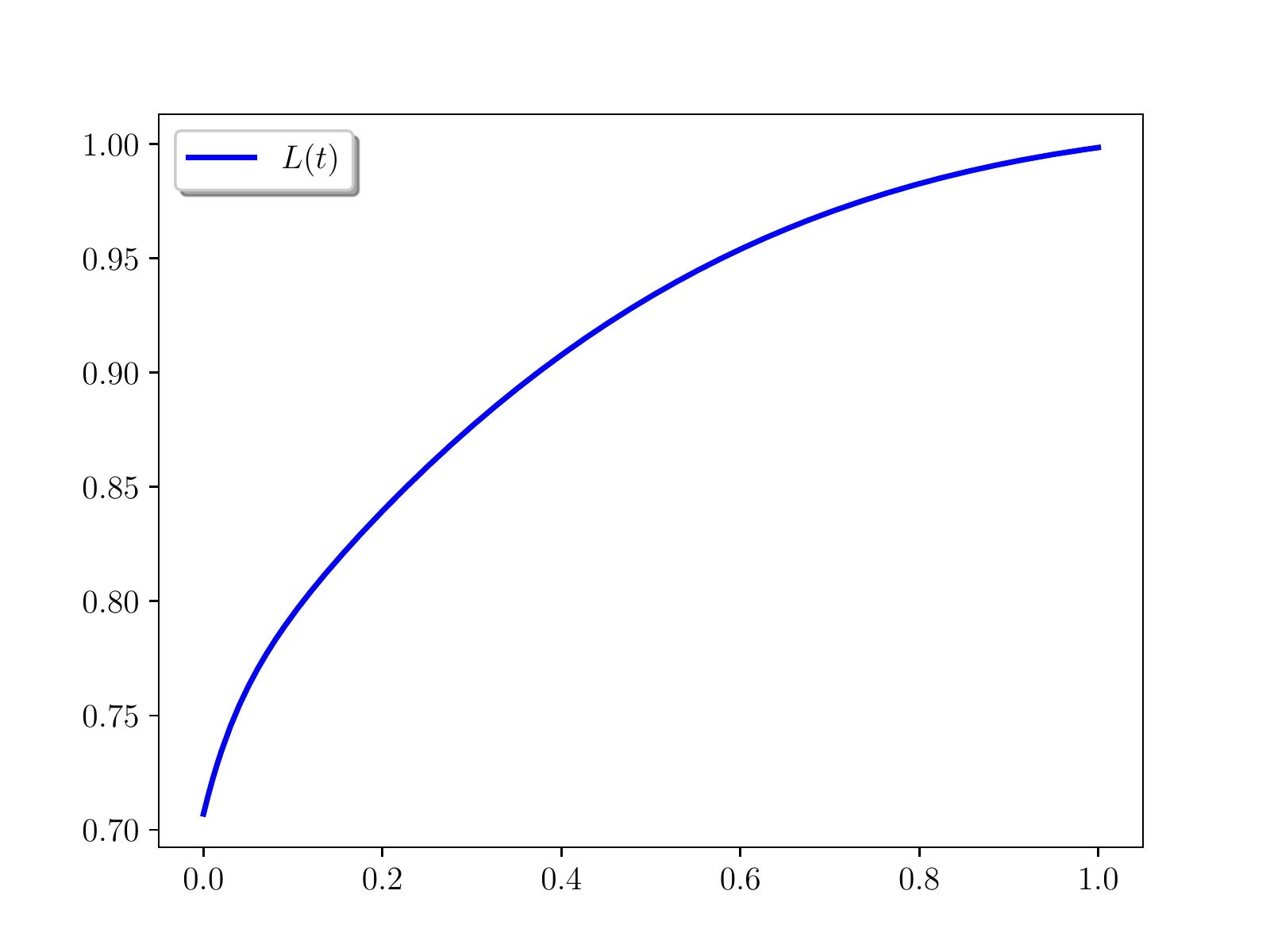} 
	\end{tabular}	
	\caption{Experiment 2. $T=1$. \textbf{(Left)} Plot of the densities in problem \eqref{s1} at different times, and \textbf{(right)}  $L(t)$.  }
	\label{exp2:state_density_length}
\end{figure}

Finally, Figure \ref{exp2:state_density_length} (left) shows the controlled densities at different times, and (right) the controlled length $L(t)$. As is the previous experiment, the results are rather accurate in what concerns the attainment of all the required initial and final conditions. 

\section{Conclusion}\label{Sec 5}
In this paper, we have analyzed the local exact controllability to the trajectories of the two components of the state variable in a one-dimensional free boundary problem governed by the Fisher-Stefan model. By applying a Neumann boundary control along the fixed boundary, we have shown that it is possible to steer the density towards that of an "ideal" trajectory and to control the position of the free boundary. Furthermore, our results also establishes the null controllability of the system. The numerical approximation has been carried out using up-to-date, deep-learning-based approaches, specifically PINNs. A comparison study between the original nonlinear problem and its linearized version has been carried out. Although the pattern of the two computed controls are similar, quantitative differences are observed, which indicates that the nonlinear Logistic term in the modeling of the problem really plays a role.     

\backmatter

%\bmhead{Supplementary information}

%If your article has accompanying supplementary file/s please state so here. 
%
%Authors reporting data from electrophoretic gels and blots should supply the full unprocessed scans for key as part of their Supplementary information. This may be requested by the editorial team/s if it is missing.
%
%Please refer to Journal-level guidance for any specific requirements.

\bmhead{Acknowledgements}

F. Periago has been supported by grant  PID2022-141957OA-C22 funded by MCIN/AEI/10.13039/501100011033, by "ERDF A way of making Europe", and the Autonomous Community of the Regi\'on of Murcia, Spain, through the programme for the development of scientific and technical research by competitive groups (21996/PI/22), included in the Regional Program for the Promotion of Scientific and Technical Research of Fundaci\'on S\'eneca -- Agencia de Ciencia y Tecnolog\'ia de la Regi\'on de Murcia.

We thank the anonymous referee for the valuable comments and suggestions.

\section*{Declarations}

%Some journals require declarations to be submitted in a standardised format. Please check the Instructions for Authors of the journal to which you are submitting to see if you need to complete this section. If yes, your manuscript must contain the following sections under the heading `Declarations':

\begin{itemize}
	\item The authors declare that they have no known competing financial interests or personal relationships that could have appeared to influence the work reported in this paper.
	\item Python scripts for the numerical experiments in the preceding section can be downloaded from 
	\url{https://github.com/fperiago/fisher_stefan_control}
%	\item Funding
%	\item Conflict of interest/Competing interests (check journal-specific guidelines for which heading to use)
%	\item Ethics approval and consent to participate
%	\item Consent for publication
%	\item Data availability 
%	\item Materials availability
%	\item Code availability 
%	\item Author contribution
\end{itemize}

\noindent


\begin{thebibliography}{00}
	
	%%%%%%%%%%%%%%%%%%%%%% R1 %%%%%%%%%%%%%%%%%%%%%%%%%%
%		\bibitem{barcena2023exact}  
%	B\'arcena-Petisco J.\,A., Fern\'andez-Cara E., Souza D.\,A. \textit{Exact controllability to the trajectories of the one-phase Stefan problem}. J. Differ. Equ. 2023;376:126--153.
	
	\bibitem{barcena2023exact}
	B\'arcena-Petisco J.A., Fern\'andez-Cara E., Souza D.A.: Exact controllability to the trajectories of the one-phase Stefan problem. J. Differ, Equ. 376, 126--153 (2023)
	
	\bibitem{behzadan2021multiplication}
	Behzadan A., Holst M.:
	Multiplication in Sobolev spaces, revisited.
	Ark. Mat. 59(2), 275--306 (2021)
	
	
	
	
 
	
	\bibitem{brosa2019extended}  
	Brosa Planella F., Please C.P., Van Gorder R.A.: Extended Stefan problem for solidification of binary alloys in a finite planar domain. SIAM J. Appl. Math. 79(3), 876--913 (2019)
	
	
	
		\bibitem{planella2021extended} 
	Brosa Planella F., Please C.P., Van Gorder R.A.: Extended Stefan problem for the solidification of binary alloys in a sphere. Eur. J. Appl. Math. 32(2), 242--279 (2021)
	
		
	
	\bibitem{canosa1973nonlinear} 
	Canosa J.: On a nonlinear diffusion equation describing population growth. IBM J. Res. Dev 17(4), 307--313 (1973)
	
	\bibitem{crank1984free} 
	Crank J.: Free and Moving Boundary Problems. Oxford University Press (1984)
	
	
		\bibitem{dalwadi2020mathematical} 
	Dalwadi M.P., Waters S.L., Byrne H.M., Hewitt I.J.: A mathematical framework for developing freezing protocols in the cryopreservation of cells. SIAM J. Appl. Math. 80(2), 657--689 (2020)
	
	\bibitem{demarque2018local}  
	Demarque R., Fern\'andez-Cara E.: Local null controllability of one-phase Stefan problems in 2D star-shaped domains. J. Evol. Equ. 18(1), 245--261 (2018)
	

	
		
	
	\bibitem{du2010spreading} 
    Du Y., Lin Z.: Spreading-vanishing dichotomy in the diffusive Logistic model with a free boundary. SIAM J. Math. Anal. 42(1), 377--405 (2010)

	\bibitem{du2011spreading} 
	Du Y., Guo Z.: Spreading--vanishing dichotomy in a diffusive Logistic model with a free boundary, II. J. Differ. Equ. 250(12), 4336--4366 (2011)
	
	\bibitem{du2012stefan} 
	Du Y., Guo Z.: The Stefan problem for the Fisher--KPP equation. J. Differ. Equ. 253(3), 996--1035 (2012)
	
    \bibitem{et2025null}
    Et-Tahri F., Chorfi S.-E., Maniar L., Boutaayamou I.:
    Null controllability of a volume--surface reaction--diffusion equation with dynamic boundary conditions.
    \emph{J. Math. Anal. Appl.} 542(2), 128793 (2025)
	
	\bibitem{fadai2020new} 
	Fadai N.T., Simpson M.J.: New travelling wave solutions of the Porous--Fisher model with a moving boundary. J. Phys. A: Math. Theor 53(9), 095601 (2020)
	
		\bibitem{fernandez2017local}  
	Fern\'andez-Cara E., de Sousa I.T.: Local null controllability of a free-boundary problem for the semilinear 1D heat equation. Bull. Braz. Math. Soc., New Ser. 48(2), 303--315 (2017)
	
	\bibitem{fernandez2013strong}  
	Fern\'andez-Cara E., M\"unch A.: Strong convergent approximations of null controls for the 1D heat equation. SeMA J. 61(1), 49--78 (2013)
	
	\bibitem{fernandez2016controllability}  
	Fern\'andez-Cara E., Limaco J., de Menezes S.B.: On the controllability of a free-boundary problem for the 1D heat equation. Syst. Control Lett. 87, 29--35 (2016)
	

	
		\bibitem{fernandez2017theoretical}  
	Fern\'andez-Cara E., Nina-Huam\'an D., N\'u\~nez-Ch\'avez M.R., Vieira F.B.: On the theoretical and numerical control of a one-dimensional nonlinear parabolic partial differential equation. J. Optim. Theory Appl. 175, 652--682 (2017)
	
	
	
	\bibitem{fernandez2019local}  
	Fern\'andez-Cara E., Hern\'andez F., L\'imaco J.: Local null controllability of a 1D Stefan problem. Bull. Braz. Math. Soc., New Ser. 50, 745--769 (2019)
	
	
	
		\bibitem{fernandez2023local}  
	Fern\'andez-Cara E., L\'imaco J., Thamsten Y., Menezes D.: Local null controllability of a quasi-linear system and related numerical experiments. ESAIM: Control Optim. Calc. Var. (2023). https://doi.org/10.1051/COCV/2023009
	
	
	
	

	
	\bibitem{fursikov1996controllability} 
	Fursikov A.V,, Imanuvilov O.Y.: Controllability of Evolution Equations. Lecture Notes Series, Vol. 34. Seoul: Seoul National University, Research Institute of Mathematics, Global Analysis Research Center (1996)
	
	\bibitem{gaffney1999modelling} 
	Gaffney E.A., Maini P.K., McCaig C.D., Zhao M., Forrester J.V.: Modelling corneal epithelial wound closure in the presence of physiological electric fields via a moving boundary formalism. Math. Med. Biol 16(4), 369--393 (1999)
	
	\bibitem{per23}  
	Garc\'ia Cervera C.J., Kessler M., Periago F.: Control of partial differential equations via physics-informed neural networks. J. Optim. Theory Appl. 196(2), 391--414 (2023)
	
	
	\bibitem{glorot2010understanding} 
	Glorot, X., and Bengio, Y.: Understanding the difficulty of training deep feedforward neural networks. In: Proceedings of the Thirteenth International Conference on Artificial Intelligence and Statistics, JMLR Workshop and Conference Proceedings pp. 249--256 (2010)
	
	\bibitem{grindrod1991patterns} 
	Grindrod P.: The Theory and Applications of Reaction-Diffusion Qquations : Patterns and waves. Clarendon Press, Oxford / Oxford Univ. Press, New York (1991)
	
	
	

	
	\bibitem{kimpton2013multiple} 
	Kimpton L.S.,  Whiteley J.P., Waters S.L., King J.R., Oliver J.M.: Multiple travelling-wave solutions in a minimal model for cell motility. Math. Med. Biol. (IMA) 30(3), 241--272 (2013)
	
	
	\bibitem{kingma2014adam} 
	Kingma, D.P., Ba, J.: Adam: A method for stochastic optimization. In: 3rd International Conference on Learning Representations (2015)
	
	\bibitem{kolmogorov1937study} 
	Kolmogorov A.N.: A study of the equation of diffusion with increase in the quantity of matter, and its application to a biological problem. Mosc. Univ. Math. Bull 1, 1--25 (1937)
	
	\bibitem{lewis2016mathematics} 
	Lewis M.A. et al.: The Mathematics Behind Biological Invasions. Springer International Publishing (2016)
	
	
	\bibitem{lions2012non}
	Lions, J. L.,  Magenes, E.: Non-Homogeneous Boundary Value Problems and Applications: Vol. 1, Vol, 181. Springer Science \& Business Media. (2012)
	
	
	
	\bibitem{mitchell2010improving} 
	Mitchell S.L, Myers T.G.: Improving the accuracy of heat balance integral methods applied to thermal problems with time dependent boundary conditions. Int. J. Heat Mass Transf. 53(17-18), 3540--3551 (2010)
	
	
	\bibitem{mitchell2009finite} 
	Mitchell S.L., Vynnycky M.: Finite-difference methods with increased accuracy and correct initialization for one-dimensional Stefan problems. Appl. Math. Comput 215(4), 1609--1621 (2009)
	
	
	
	
	
	
	\bibitem{murray2002mathematical} 
	Murray J.D. Mathematical biology: I. An Introduction. Interdisciplinary Applied Mathematics. 17 (2002)
	
	\bibitem{huaman2023local}  
	Nina Huam\'an D., N\'u\~nez-Ch\'avez M.R., L\'imaco J., Carvalho P.P.: Local null controllability for the thermistor problem. Nonlinear Anal. 236, 113330 (2023)
	
	
	\bibitem{PENG2022106067}  
	Peng W.Q., Pu J.C., Chen Y.: PINN deep learning method for the Chen--Lee--Liu equation: Rogue wave on the periodic background. Commun. Nonlinear Sci. Numer. Simul. 105, 106067 (2022)
	
	
		\bibitem{pinkus1999approximation} 
	Pinkus A.: Approximation theory of the MLP model in neural networks. Acta Numer. 8, 143--195 (1999)
	
	\bibitem{rai19} 
	Raisi M., Perdikaris P., Karniadakis G.E.: Physics-informed neural networks: A deep learning framework for solving forward and inverse problems involving nonlinear partial differential equations. J. Comput. Phys 378, 686--707 (2019)
	
	
	
	\bibitem{shuttleworth2019multiscale} 
	Shuttleworth R., Trucu D.: Multiscale modelling of fibres dynamics and cell adhesion within moving boundary cancer invasion. Bull. Math. Biol 81, 2176--2219 (2019)
	
	\bibitem{trelat2005controle}  
	Tr\'elat E.: Contr\^ole Optimal : Th\'eorie \& Applications. Vuibert, Paris (2005)
	
	\bibitem{ward1997mathematical} 
	Ward J.P., King J.R.: Mathematical modelling of avascular-tumour growth. Math. Med. Biol. (IMA) 14(1), 39--69 (1997)
	
	
\bibitem{ZHANG2024108229}  
Zhang Q., Qiu C., Hou J., Yan W.: Advanced Physics-informed neural networks for numerical approximation of the coupled Schr\"odinger--KdV equation. Commun. Nonlinear Sci. Numer. Simul 138, 108229 (2024)

	
	
	
\end{thebibliography}
\end{document}